\documentclass{amsart}
\usepackage[english]{babel}
\usepackage{amsaddr}
\usepackage[utf8]{inputenc}
\usepackage{enumitem}
\usepackage{esint}

\usepackage[backend=biber]{biblatex}
\usepackage{lmodern}

\usepackage{amsmath}

\usepackage{amssymb}

\usepackage{bbm}
\usepackage{csquotes}
\usepackage{fancyhdr}

\numberwithin{equation}{section}

\usepackage{xcolor}
\usepackage{geometry}
\newtheoremstyle{theorem2}
{8pt}
{8pt}
{\itshape}
{}
{\bfseries}
{.}
{.5em}
{}

\newtheoremstyle{definition2}
{8pt}
{8pt}
{\itshape}
{}
{\bfseries}
{.}
{.5em}
{}

\theoremstyle{theorem2}

\newtheorem{lemma}{Lemma}[section]
\newtheorem{teo}[lemma]{Theorem}

\newtheorem{cor}[lemma]{Corollary}
\newtheorem{prop}[lemma]{Proposition}

\theoremstyle{definition2}

\newtheorem{deff}[lemma]{Definition}
\newtheorem{Remark}[lemma]{Remark}
\newtheorem{exa}[lemma]{Example}

\begin{document}
	
\title{A continuous model of transportation\\ in the Heisenberg group}

\author{Michele Circelli}
\address{University of Bologna, Department of Mathematics, 40126-Bologna (Italy)}
\email{michele.circelli2@unibo.it}

\author{Albert Clop}
\address{University of Barcelona, Department of Mathematics and Computer Science, 08007-Barcelona (Catalonia)}
\email{albert.clop@ub.edu}

\begin{abstract}
We present a minimization problem with a horizontal divergence-type constraint in the Heisenberg group. Our study explores its dual formulation and examines its relationship with the congested optimal transport problem, for $1 < p < +\infty$, as well as the Monge-Kantorovich problem, in the limite case $p=1$.
\end{abstract}

\maketitle

\section{Introduction}\label{introduction}

The optimal transport problem, originally posed by G. Monge in \cite{monge1781memoire} in the late 18th century, has become a highly active research area due to its numerous applications, as highlighted in \cite{evans1997partial} and \cite{rachev2006mass}. In the Monge-Kantorovich formulation, which dates back to 1942, we are given two probability measures $\mu$ and $\nu$ on a metric space $(X,d)$ and aim to solve the following minimization problem:
\begin{equation}\label{introd22gennaio}
    \inf_{\gamma\in\Pi(\mu,\nu)}\int_{X\times X} c(x,y)d\gamma(x,y),
\end{equation}
where the set $\Pi(\mu,\nu)$ consists of \textit{transport plans}, i.e. probability measures $\gamma$ on the product space $X\times X$, with prescribed marginals $\mu$ and $\nu$. The cost function $c\geq0$ is typically a power of the distance $d$. Solutions exist under quite general assumptions, such as when $X$ is a Polish space and $c$ is lower semicontinuous and bounded from below.

Initially, the existence of solutions to \eqref{introd22gennaio} induced by measurable maps was extensively studied within Euclidean and Riemannian settings by many researchers (see, for instance, \cite{ambrosio2021lectures}, \cite{Santambrogiolibro}, \cite{Feldman} and references therein for comprehensive overviews). Later, the research focus shifted towards sub-Riemannian frameworks, especially Carnot groups, which arise as limits of Riemannian manifolds in a very precise sense (see \cite{Gromov}).

The Heisenberg group $\mathbb{H}^n$ is the simplest model of non-commutative Carnot group: it is characterized by the selection of $2n$ vector fields, $X_1, \ldots, X_{2n}$, in $\mathbb{R}^{2n+1}$, known as horizontal vector fields, and a metric on the sub-bundle $H\mathbb{H}^n$ they generate on the tangent bundle $T\mathbb{H}^n$. Since the number of vector fields is less than the dimension of the space, the metric is non-Riemannian at every point. Moreover, displacement occurs only along integral curves of horizontal vector fields, also called horizontal curves. The presence of this metric naturally leads to the definition of an intrinsic sub-Riemannian distance $d_{SR}$.

The first significant work on the Monge-Kantorovich problem in $\mathbb{H}^n$
\begin{equation}\label{introd22gennaio2}
    \inf_{\gamma\in\Pi(\mu,\nu)}\int_{\mathbb{H}^n\times\mathbb{H}^n}d_{SR}(x,y)^\alpha d\gamma(x,y)
\end{equation}
is \cite{Rigot}, where L. Ambrosio and S. Rigot proved the existence and the uniqueness of the solution to \eqref{introd22gennaio2} for $\alpha=2$, and that this solution is induced by a map; subsequently, this result was extended to general sub-Riemannian manifolds in \cite{agrachevlee}, \cite{FigalliRifford}, and to essentially non-branching metric measure spaces in \cite{Cavalletti2}. The first work on the problem \eqref{introd22gennaio2} for $\alpha=1$ is \cite{DePascale2}, where the authors proved the existence of solutions induced by maps. Subsequently, the analogous result was proven in more general measure metric spaces, see for instance \cite{Cavalletti3} and \cite{cavalletti2018overview}. One of the main obstacles in transitioning from the Euclidean to the sub-Riemannian setting is the loss of uniqueness of geodesics, even locally, and hence, addressing the problem in this new framework required new and very different approaches due to the unique geometry of the space.

In this direction we continue, in this paper, the study initiated in \cite{circelli2023transport} on the congested optimal transport problem in $\mathbb{H}^n$: this study is motivated by the proposal of Petitot and Tondut in \cite{Petitot} to model the functional geometry of the visual cortex (see also \cite{cittisarti}, where this model was expressed using sub-Riemannian geometry). The visual signal propagates along the connectivity, leading to congestion phenomena that can be modeled with the congested optimal transport problem in $\mathbb{H}^n$ (see \cite[Introduction]{circelli2023transport} and references therein). 

The congested optimal transport problem was first introduced in the Euclidean setting by Carlier et al. in \cite{Santambrogio1}. In the Heisenberg group it reads as follows: given two probability measures $\mu$ and $\nu$ over the closure of a bounded domain $\Omega$, and a cost function $G:\Omega\times\mathbb{R}_+\rightarrow\mathbb{R}_+$, which is convex and exhibits $p$-growth in the second variable, it consists of solving 
\begin{equation}\label{intro3}
	\min_{Q\in\mathcal{Q}_H^p(\mu,\nu)}\int_\Omega G(x,i_Q(x))dx.
\end{equation}
Here, the set $\mathcal{Q}_H^p(\mu,\nu)$ consists of probability measures $Q$ over the space $H$ of horizontal curves, with prescribed projections $\mu$ and $\nu$ and for which the associated \textit{horizontal traffic intensity} $i_Q$
$$\int_{\overline\Omega} \varphi(x) d i_{Q}(x):=\int_H\left(\int_0^1\varphi(\sigma(t))|\dot{\sigma}(t)|_Hdt\right) d Q(\sigma),\quad \forall \varphi \in C(\overline\Omega)$$
is a $p$-summable function; $|\cdot|_H$ denotes a norm on the sub-bundle $H\mathbb{H}^n$. In \cite{circelli2023transport} the authors proved that solutions to \eqref{intro3} are \textit{equilibrium configurations of Wardrop type}. Although \eqref{intro3} is a convex minimization problem, and hence the existence of solutions follows from the direct method in the calculus of variations, it is also a minimization problem over a set of measures over curves, and hence it has two layers of infinite dimensions. 
Therefore, in the spirit of \cite{brasco2013congested} and \cite{santambrogio2014dacorogna}, we investigate the relation between \eqref{intro3} and the following more tractable problem:
\begin{equation}\label{intro1}
	\min_{\textnormal{\textbf{w}}}\int_\Omega \mathcal{G}(x,\textnormal{\textbf{w}}(x))dx,
\end{equation}
where $\mathcal{G}(x,\textnormal{\textbf{w}})=G(x,|\textnormal{\textbf{w}}|_H)$,
and the minimum is taken among all $p$-summable horizontal vector fields $\textnormal{\textbf{w}}:\Omega\rightarrow H\mathbb{H}^n$, with prescribed horizontal divergence
\begin{equation}\label{intro4}
	\textnormal{div}_H\textnormal{\textbf{w}}=\mu-\nu.
\end{equation}
This problem is the Heisenberg version of the well-known \textit{continuous model of transportation} introduced in the Euclidean setting by Beckmann in \cite{Beckmann2}. Following the approach \cite{brasco2013congested} and \cite{santambrogio2014dacorogna}, which relies on the Dacorogna-Moser flow argument, together with a standard regularization procedure, we prove that the two problems \eqref{intro1} and \eqref{intro3} are actually equivalent. Moreover, in some sense, minimizers to \eqref{intro3} and \eqref{intro1} correspond to each other, see Proposition \ref{21marzo251} and Theorem \ref{9dicembre2}. The main difference here with the Euclidean case is that, in order to apply the Dacorogna-Moser flow argument, we must pass through a Riemannian approximation of the Heisenberg group.

Using standard duality arguments, \eqref{intro1} can be formulated as a classical problem in the calculus of variations, and hence it is strongly related to the regularity theory for PDEs in the Heisenberg group. More precisely, in the spirit of \cite{santambrogioregularity}, we show in Theorem \ref{dualteo} that \eqref{intro1} is equivalent to the following problem:
\begin{equation}\label{intro2}
	\max _{\varphi\in HW^{1,q}(\Omega)}\left\{-\int_\Omega\varphi d(\mu-\nu)-\int_\Omega  \mathcal{G}^*(x,\nabla_H\varphi(x))dx\right\},
\end{equation}
where $\mathcal{G}^*$ denotes the Legendre transform of $\mathcal{G}$ in the second variable.

Here, $\nabla_H\varphi$ denotes the vector of derivatives of $\varphi$ in horizontal directions, and $HW^{1,q}(\Omega)$ is its associated horizontal Sobolev space, where the exponent $q$ is the dual of $p$ (see Section 2 for its definition). Moreover, if $(\textnormal{\textbf{w}}_0,\varphi_0)$ is a pair of minimizer for \eqref{intro1} and maximizer for \eqref{intro2}, then, if $\mathcal{G}^*$ is differentiable, we have
\begin{equation*}
	\textnormal{\textbf{w}}_0(x)= D\mathcal{G}^*(x,\nabla_H\varphi_0(x)) \text{ a.e. in }\Omega,
\end{equation*}
where $D\mathcal{G}^*$ is the Euclidean gradient of $\mathcal{G}^*$ in $\mathbb{R}^{2n+1}$, and $\varphi_0$ solves the associated Euler-Lagrange equation of horizontal divergence-type:
\begin{equation}\label{intro6}
    \textnormal{div}_H \left(D\mathcal{G}^*(x,\nabla_H\varphi(x))\right)=\mu-\nu.
\end{equation}

In the Euclidean setting, the equivalence between \eqref{intro3}, \eqref{intro1} and \eqref{intro2} was first investigated in \cite{Brasco} for the autonomous case. More recently, this equivalence was extended to the non-autonomous case in \cite{Brasco2}, where the authors also assume that the right-hand side in \eqref{intro4} belongs to the dual of some Sobolev space, whose elements are not measures in general.

In the limiting case $p=1$, the simplest version of the problem \eqref{intro1} reads as:
\begin{equation}\label{L1beckmann}
	\inf_{\textnormal{\textbf{w}}\in L^1 }\left\{\int_\Omega|{\textbf{w}}(x)|_Hdx:\textnormal{div}_H\textnormal{\textbf{w}}=\mu-\nu\right\}.
\end{equation}
As observed in \cite[Subsection 4.2.2]{Santambrogiolibro}, this problem is not well-posed a priori: due to the non-reflexivity of $L^1$, there may not exist any $L^1$ horizontal vector field minimizing the $L^1$ norm under a horizontal divergence constraint. This is why one states the problem \eqref{L1beckmann} in the larger space of compactly supported horizontal vector measures, whose horizontal divergence is the signed Radon measure $\mu-\nu$:
\begin{equation}\label{intro5}
	\min\left\{\|\textnormal{\textbf{w}}\|:\textnormal{div}_H \textnormal{\textbf{w}}=\mu-\nu\right\},
\end{equation}
where $\|\textnormal{\textbf{w}}\|$ is the total variation of the vector measure $\textnormal{\textbf{w}}$. In the spirit of \cite[Section 4]{Santambrogiolibro}, we show in Theorem \ref{7dicembre1} that \eqref{intro5} equals the value of the Monge-Kantorovich problem \eqref{introd22gennaio2}, with $\alpha=1$. A key role in the proof of this equivalence is played by the transport density, which measures how much transport is taking place in any region of $\mathbb{H}^n$. We delve deeper into the study of the properties of the transport density in $\mathbb{H}^n$, initiated in \cite{circelli2023transport}, explicitly describing its relation with Kantorovich potentials. To do so, we adapt some results on the differentiability of the Kantorovich potentials present in \cite{Feldman}, which were specific to the Riemannian case. While adapting these results, we strongly rely on results from \cite{Rigot} and \cite{DePascale1}, which were based on the explicit computation of geodesics, the smoothness of the pointed sub-Riemannian distance, and the measure contraction property.

As byproducts of the previous equivalence, we obtain two important facts: first, if either $\mu$ or $\nu$ is absolutely continuous with respect to the Lebesgue measure (which is a Haar measure for the group structure), then \eqref{L1beckmann} admits a solution. Second, one can conclude that the problem \eqref{intro5} also admits a dual formulation (as it does for its congested version \eqref{intro1}), which is given by the well-known \textit{Kantorovich duality theorem}:
\begin{equation*}
\max_{u \in HW^{1,\infty}} \left\{ \int u \, d(\mu - \nu) : \|\nabla_H u\|_{\infty} \leq 1 \right\}.
\end{equation*}

In addition, in Theorem \ref{dynteo} we prove that the Lagrangian reformulation of \eqref{intro5} is:
\begin{equation*}
	\min_{Q}\int_H\ell_{SR}(\sigma)dQ(\sigma)=\min_Q\int di_Q(x),
\end{equation*}
where $\ell_{SR}(\sigma)$ denotes the \textit{sub-Riemannian length} (see Section 2 for its definition). We emphasize that in the proof, the geometry of the space plays no role, so the result holds in the more general setting of Polish metric spaces.

The paper is organized as follows: In Section \ref{Preliminaries}, we introduce the Heisenberg group and recall some well-known results on optimal transport theory in this setting. In Section \ref{Transportrays}, we prove some results about the transport set and the Kantorovich potential, which will be useful later. Section \ref{Beckmannp1} deals with the aforementioned three problems in the limit case $p=1$, while the last section addresses the equivalence of the three problems for $1<p<+\infty$. First, we introduce both \eqref{intro1} and \eqref{intro2} in their natural functional analytic setting; second, we prove the equivalence between \eqref{intro1} and the congested optimal transport problem \eqref{intro3}.

\section{Preliminaries}\label{Preliminaries}

\subsection{The Heisenberg group $\mathbb{H}^n$}
Let $n\geq1$. The $n$-th \textit{Heisenberg group} $\mathbb{H}^n$ is the connected and simply connected Lie group, whose Lie algebra $\mathfrak{h}^n$ is stratified of step 2: i.e. it is the direct sum of two linear subspaces
$$ \mathfrak{h}^n= \mathfrak{h}^n_1 \oplus \mathfrak{h}^n_2,$$

where $\mathfrak{h}_1= \mathrm{span} \{X_1, \dots, X_n, X_{n+1}, \dots, X_{2n}\}$ is the horizontal layer, $\mathfrak{h}_2=\mathrm{span} \{X_{2n+1} \}$ and the only non-trivial bracket-relation between the vector fields $X_1, \dots, X_n, X_{n+1}, \dots, X_{2n}, X_{2n+1}$ is 
\begin{equation*}
	[X_j, X_{n+j}]=X_{2n+1},\quad \forall j=1,\ldots n.
\end{equation*}
The homogeneous dimension of $\mathbb{H}^n$ is 
\begin{equation*}\label{homogdim}
	N:=\sum_{i=1}^2i\dim(\mathfrak{h}^n_i)=2n+2.
\end{equation*}
The exponential map $\mathrm{exp}: \mathfrak{h}^n\to \mathbb{H}^n$ is a diffeomorphism, hence it induces a global diffeomorphism between $\mathbb{H}^n$ and $\mathbb{R}^{2n+1}$:
\begin{equation*}
	\mathbb{H}^n\ni\exp\left(x_1X_1+x_2X_2+\ldots+x_{2n}X_{2n}+x_{2n+1}X_{2n+1}\right)\longleftrightarrow\left(x_1,\ldots,x_{2n+1}\right)\in\mathbb{R}^{2n+1}.
\end{equation*} 
In this system of coordinates the group law will be given, through the Baker-Campbell-Hausdorff formula which defines a group structure on $\mathfrak{h}^n$, by  
\begin{equation*}\label{grouplaw}
	x \cdot y := \bigg{(}x_1+y_1,\ldots,x_{2n}+y_{2n}, x_{2n+1}+y_{2n+1}+\frac{1}{2}\sum_{j=1}^{2n} (x_jy_{n+j}-x_{n+j}y_j)\bigg{)},
\end{equation*}
the unit element $e\in\mathbb{H}^n$ is $0_{\mathbb{R}^{2n+1}}$, the \textit{center} of the group is 
\begin{equation*}
	L := \{(0,\ldots,0,x_{2n+1}) \in \mathbb{H}^n ;\; x_{2n+1} \in \mathbb{R}\}
\end{equation*}
and the vector fields $X_1,\ldots,X_{2n+1}$, left invariant w.r.t. \eqref{grouplaw}, in coordinates read as
\begin{equation*}
	\begin{cases}
		X_j := \partial_{x_j} -\frac{x_{n+j}}{2} \partial_{x_{2n+1}} \,, j=1,\dots,n,\\
		X_{n+j} := \partial_{x_{n+j}}+\frac{x_j}{2} \partial_{x_{2n+1}},\ \ j=1,\dots,n,\\
		X_{2n+1}=\partial_{x_{2n+1}}.	
	\end{cases}
\end{equation*}
According to the two step-stratification of the algebra,  $\mathbb{H}^n$ is naturally endowed with a family of intrinsic non-isotropic dilations that reads as 
\begin{equation*}
	\delta_\tau((x_1,\ldots,x_{2n+1})):= (\tau x_1,\ldots,\tau x_{2n},\tau^2x_{2n+1}),\quad \forall x\in\mathbb{H}^n, \forall \tau>0.
\end{equation*}
See \cite[Section 2]{Bonfiglioli} for more details about this topic.

The horizontal layer $\mathfrak{h}^n_1$ of the algebra defines a sub-bundle $H\mathbb{H}^n$ of the tangent bundle $T\mathbb{H}^n$, whose fibre at a point $x \in \mathbb{H}^n$ is
$$H_x \mathbb{H}^n=\mathrm{span} \{ X_1(x), \dots, X_n(x), X_{n+1}(x), \dots, X_{2n}(x) \}.$$

One can equip $\mathbb{H}^n$ with a left-invariant sub-Riemannian metric as follows: we fix an inner product $\langle \cdot, \cdot \rangle_H$ on $\mathfrak{h}_1^n$ that makes $\{X_1, \dots, X_n, X_{n+1}, \dots, X_{2n}\}$ an orthonormal basis and we denote by $|\cdot|_H$ the norm associated with such a scalar product. We keep the same notation both for the corresponding scalar product and norm on each fibre $H_x\mathbb{H}^n$, $x\in\mathbb{H}^n$.

An absolutely continuous curve $\sigma\in AC([0,1],\mathbb{R}^{2n+1})$ is said to be \textit{horizontal} if its velocity vector $\dot{\sigma}(t)$ belongs to $H_{\sigma(t)}\mathbb{H}^n$ at almost every $t\in[0,1]$. Its \textit{sub-Riemannian length} is 
\begin{equation}\label{horlength}
	\ell_{SR}(\sigma):=\int_0^1|\dot{\sigma}(t)|_H dt.
\end{equation}
Let us just remark that $\ell_{SR}$ is invariant under reparametrizations of $\sigma$.

The Rashevsky-Chow's theorem guarantees that the \textit{sub-Riemannian distance} 
\begin{equation*}\label{CCdistance}
	d_{SR}(x,y):= \inf \  \left\{  \ell_{SR}(\sigma):\sigma \text{ is horizontal}, \ \sigma(0)=x, \ \sigma(1)=y \right\},\quad \forall x,y\in\mathbb{H}^n
\end{equation*}
is well-defined and it induces the Euclidean topology, see \cite[Theorem 3.31]{agrachev2019comprehensive}. In particular, $(\mathbb{H}^n,d_{SR})$ is a polish space. 

If $x\in\mathbb{H}^n$ and $r>0$, we denote by
\begin{equation*}
	B(x,r):=\left\{y\in\mathbb{H}^n: d_{SR}(x,y)< r\right\}
\end{equation*}
the \textit{sub-Riemannian ball} centered in $x$ with radius $r$.
The Lebesgue measure $\mathcal{L}^{2n+1}$ is a Haar measure on $\mathbb{H}^n\simeq\mathbb{R}^{2n+1}$; we will often denote it by $dx$.

For any open set $\Omega\subseteq\mathbb{H}^n$ and any measurable function $\varphi:\Omega\to\mathbb{R}$, we denote by 
$$\nabla_H\varphi=\sum_{i=1}^{2n}X_i\varphi X_i$$ its \textit{horizontal gradient}, in the sense of distributions. 

We say that a measurable map $\varphi:\Omega \rightarrow \mathbb{R}$ is \textit{Pansu differentiable} at $x \in \mathbb{H}^n$ if there exists an homogeneous group homomorphism $L :\mathbb{H}^n \rightarrow \mathbb{R}$ such that
\begin{equation*}
	\lim_{y \rightarrow x} \frac{\varphi(y) - \varphi(x) - L(x^{-1} \cdot y)}{d_{SR}(x,y)} = 0,
\end{equation*}
where a group homomorphism $L:\mathbb{H}^n\to\mathbb{R}$ is \textit{homogeneous} if $L(\delta_\tau(x))=\tau L(x)$, for any $x\in\mathbb{H}^n$ and any $\tau>0$.

If the map $L$ exists, it is unique and will be denoted by $D_H \varphi(x)$.

If $\varphi:\Omega\rightarrow \mathbb{R}$ is Pansu differentiable at $x \in \Omega$ then $\varphi$ has directional derivatives at $x$ along all the directions $X_j,\forall j=1,\ldots,2n$. Moreover, if $y=(y_1,\dots, y_{2n+1})\in \mathbb{H}^n$ then 
\begin{equation}\label{9maggio}
	D_H\varphi(x)(y)=\sum_{j=1}^{2n} y_j X_j\varphi(x)=\left\langle\nabla_H\varphi(x),\pi_{x}(y)\right\rangle_H ,
\end{equation}
where the map $y\mapsto\pi_{x}(y)$ is the smooth section of $H\mathbb{H}^n$ defined as
\begin{equation}
	\pi_{x}(y):=\sum_{i=1}^{2n}y_iX_i(y),
\end{equation}
see \cite[Proposition 5.6]{Franchi1}. If $\varphi\in C^\infty$ and $\sigma$ is an horizontal curve, a simple computation based on the classical chain rule implies that
\begin{equation}\label{14marzo}
	\frac{d}{dt}\varphi(\sigma(t))\bigg{|}_{t=s}=\langle\nabla_H\varphi(\sigma(s)),\dot\sigma(s)\rangle_H,
\end{equation}
for a.e. $s\in[0,1]$.

For every $1\leq q\leq\infty$, the space
\begin{equation}\label{horsob}
	HW^{1,q}(\Omega):=\left\{\varphi:\Omega\to\mathbb{R} \text{ measurable : }\varphi\in L^{q}(\Omega),\nabla_H\varphi\in L^{q}(\Omega, H\Omega)\right\},
\end{equation}
is a Banach space equipped with the norm
\begin{equation*}
	\|\varphi\|_{HW^{1,q}(\Omega)}:=\|\varphi\|_{L^q(\Omega)}+\|\nabla_H \varphi\|_{L^q(\Omega,H\Omega)}.
\end{equation*}

\subsubsection{Geodesics in $\mathbb{H}^n$}

We call \textit{minimizing horizontal curve} any curve $\sigma\in AC([0,1],\mathbb{R}^{2n+1})$, which is horizontal and such that
\begin{equation*}
	\ell_{SR}(\sigma)=d_{SR}(\sigma(0),\sigma(1)),
\end{equation*}
and \textit{geodesic} any minimizing horizontal curve parametrized proportionally to the arc-length.

The space $(\mathbb{H}^n, d_{SR})$ is a complete and locally compact metric length space, therefore it is geodesic \cite[Theorem 2.5.28]{burago2022course}.

Usually, in this kind of spaces, geodesics between two points may not be  unique: in the case of the Heisenberg group geodesics can be computed explicitly and it is possible to detect a set in which these are unique.

We denote by
\begin{equation} \label{KAPPA}
	E:= \{(x,y)\in \mathbb{H}^n\times \mathbb{H}^n;\,\, x^{-1}\cdot y \not \in L\},
\end{equation}
then it holds the following characterization for geodesics parametrized on $[0,1]$.

\begin{teo}\label{geod} 
A non trivial geodesics parametrized on $[0,1]$ and starting from $0$ is the restriction to $[0,1]$ of the curve $\sigma_{\chi,\theta}(t)=\left(x_1(t),\ldots,x_{2n+1}(t)\right)$ either of the form
\begin{align}\label{geodform}
	&x_j(t)=\frac{\chi_j\sin(\theta s)-\chi_{n+j}\left(1-\cos(\theta s)\right)}{\theta},\quad j=1,\ldots,n\\
    &x_{n+j}(t)=\frac{\chi_{n+j}\sin(\theta s)+\chi_{j}\left(1-\cos(\theta s)\right)}{\theta},\quad j=1,\ldots,n\\
    &x_{2n+1}(t)=\frac{|\chi|^2}{2\theta^2}\left(\theta s-\sin(\theta s)\right),
\end{align}
for some $\chi \in \mathbb{R}^{2n}\setminus\{0\}$ and $\theta\in [-2\pi,2\pi]\setminus\left\{0\right\}$, or of the form
\begin{equation*}
	\left(x_1(t),\ldots,x_{2n+1}(t)\right)=\left(\chi_1t,\ldots,\chi_{2n}t,0\right), 
\end{equation*}
for some $\chi \in \mathbb{R}^{2n}\setminus\{0\}$ and $\theta=0$. In particular, it holds $$|\chi|_{\mathbb{R}^{2n}}=|\dot\sigma_{\chi,\theta}|_H=d_{SR}(x,\sigma_{\chi,\theta}(1)).$$
	
Moreover it holds:
\begin{enumerate}
    \item For all $(x,y)\in E$, there is a unique geodesic $x\cdot\sigma_{\chi,\theta}$ parametrized on $[0,1]$ between $x$ and $y$, for some $\chi \in \mathbb{R}^{2n}\setminus\{0\}$ and some $\theta\in(-2\pi,2\pi)$.
    \item If $(x,y)\not\in E$, then $x^{-1}\cdot y = (0,\ldots,0,z_{2n+1})$ for some $z_{2n+1}\in \mathbb{R}\setminus\{0\}$, there are infinitely many geodesics parametrized on $[0,1]$ between $x$ and $y$. These curves are all curves of the form $x\cdot \sigma_{\chi,2\pi}$, if $z_{2n+1}>0$,  $x\cdot \sigma_{\chi,-2\pi}$, if $z_{2n+1}<0$, for all $\chi \in \mathbb{R}^{2n}$ such that $|\chi| = \sqrt{4\pi|z_{2n+1}|} $.
\end{enumerate}
\end{teo}

From the description of geodesics given in Theorem \ref{geod} it follows that $(\mathbb{H}^n, d_{SR})$ is a non-branching metric space: any two geodesics which coincide on a non-trivial interval coincide on the whole intersection of their intervals of definition.

We denote by \textnormal{Geo}$(\mathbb{H}^n)$ the space of geodesics parametrized on $[0,1]$:
\begin{equation*}
    \textnormal{Geo}(\mathbb{H}^n):=\left\{\sigma\in AC\left([0,1],\mathbb{R}^{2n+1}\right):\sigma\textnormal{ is a geodesic}\right\}.
\end{equation*}

Let us just remark that there exists a map 
\begin{equation}\label{19marzo1}
	S:\mathbb{H}^n\times\mathbb{H}^n\rightarrow\textnormal{Geo}(\mathbb{H}^n),
\end{equation}
which is $\gamma$- measurable for any positive Borel measure $\gamma$ on $\mathbb{H}^n\times\mathbb{H}^n$, and that associates with any pair of points $(x,y)\in \mathbb{H}^n\times\mathbb{H}^n$ a geodesic $S(x,y):=\sigma_{x,y}\in\textnormal{Geo}(\mathbb{H}^n)$, between $x$ and $y$. The existence of such a map follows from the theory of Souslin sets and general theorems about measurable selections, see for instance \cite[Theorem 6.9.2 and Theroem 7.4.1]{Bogachev}. Moreover $S$ is jointly continuous on $E$. If $e_t$ is the evaluation map at time $t\in[0,1]$, then the map
\begin{equation*}
	S_t:=e_t\circ S:\mathbb{H}^n\times\mathbb{H}^n\rightarrow\mathbb{H}^n
\end{equation*}
associates to any two points $x,y\in\mathbb{H}^n$ the point $S_t(x,y):=\sigma_{x,y}(t)$ of $\mathbb{H}^n$ at distance $t \, d_{SR}(x,y)$ from $x$ on the selected geodesic $\sigma_{x,y}$ between $x$ and $y$.

\subsubsection{$H\Omega$-valued finite Radon measures}\label{vectormeasuressubsection}

Let $\Omega\subseteq \mathbb{H}^n$ be a set, and denote by $H\Omega$ the restriction of the horizontal bundle $H\mathbb{H}^n$ to the set $\Omega$. We denote by  $C_c(\Omega, H\Omega)$ the class of continuous horizontal vector fields with  compact support in $\Omega$, and by $C_0(\Omega, H\Omega)$ its completion with respect to the uniform norm:  $$\|\phi\|_\infty=\sup\{|\phi(x)|_H\,:\, x\in \Omega\},\quad \phi\in C(\Omega,H\Omega).$$
The normed space $(C_0(\Omega, H\Omega),\|\cdot\|_\infty)$ is a Banach space.

A $H\Omega$-\textit{valued finite Radon measure} $\textbf{w}$ \textit{on} $\Omega$ is a continuous linear functional on $C_0(\Omega,H\Omega)$
\begin{equation*}
    \textbf{w}\in C_0(\Omega,H\Omega)',
\end{equation*}
and its action against a horizontal vector field $\phi\in C_0(\Omega,H\Omega)$ is denoted by
\begin{equation*}
    \int_\Omega\phi\cdot d\textbf{w}.
\end{equation*}
We denote by
\begin{equation*}
    \mathcal{M}(\Omega, H\Omega):=C_0(\Omega,H\Omega)'
\end{equation*}
the space of all $H\Omega$-valued finite Radon measures on $\Omega$ and we equip it with the dual norm
\begin{equation*}
    \|\textbf{w}\|_{\mathcal{M}(\Omega,H\Omega)}:=\|\textbf{w}\|_{\left(C_0(\Omega,H\Omega)\right)'}=\sup\left\{\int_\Omega\phi\cdot d\textbf{w}:\phi\in C_0(\Omega,H\Omega),\|\phi\|_\infty\leq1 \right\}.
\end{equation*}
Since the horizontal bundle $H\Omega$ has a global trivialization, it follows from the Riesz theorem (see for instance \cite{AmbFusPal}) that $\mathcal{M}(\Omega,H\Omega)$ can be identified with the space $\mathcal{M}(\Omega)^{2n}$ of $\mathbb{R}^{2n}$-valued finite Radon measure through the identification
$$\int_{\Omega} \phi\cdot d\textbf{w}=\sum_{i=1}^{2n}\int_{\Omega}\phi_i d\textbf{w}_i,
\hspace{1cm}\forall \phi=\sum_{i=1}^{2n}\phi_iX_i\in C_0(\Omega,H\Omega),$$
for some unique 
$(\textbf{w}_1,\dots,\textbf{w}_{2n})\in\mathcal{M}(\Omega)^{2n}$. Moreover, it holds
$$
\|\textbf{w}\|_{\mathcal{M}\left(\Omega,H\Omega\right)}=|\textbf{w}|(\Omega),
$$
where $|\textbf{w}|\in\mathcal{M}_+(\Omega)$ denotes the \textit{total variation} of the vector measure $(\textbf{w}_1,\dots,\textbf{w}_{2n})$.

Finally, if $\Omega$ is compact, then $C_c(\Omega,H\Omega)=C_0(\Omega,H\Omega)=C(\Omega,H\Omega)$ and 
\begin{equation*}
	\mathcal{M}(\Omega, H\Omega)=\left(C(\Omega, H\Omega)\right)'.
\end{equation*}

\subsubsection{Mollification in $\mathbb{H}^n$}
Let us consider a mollifier for the group structure, i.e. a function $\rho\in C_c^{\infty}(\mathbb{H}^n)$, such that $\rho\geq0$, $\textnormal{supp}(\rho)\subseteq B(0,1)$ and $\int_{\mathbb{H}^n}\rho(x)dx=1$; for any $\varepsilon>0$ we denote by
\begin{equation*}
	\rho_\varepsilon(x):=\varepsilon^{-N}\rho\left(\delta_{1/\varepsilon}(x)\right).
\end{equation*}
If $\varphi\in L^1_{loc}$, then 
\begin{equation*}
	\rho_\varepsilon\ast\varphi(x):=\int_{\mathbb{H}^n}\rho_\varepsilon(x\cdot y^{-1})\varphi(y)dy,
\end{equation*}
is smooth and it enjoys many standard properties, see for instance \cite{Stein1}. 

More generally, one may also mollify Radon measures: given a finite Radon measure $\lambda\in\mathcal{M}(\mathbb{H}^n)$, resp. a finite vector Radon measure $\textbf{w}=\left(\textbf{w}_1,\ldots,\textbf{w}_{2n}\right)\in\mathcal{M}(\mathbb{H}^n,H\mathbb{H}^n)$, then the mollified functions $\lambda^\varepsilon\in C^{\infty}(\mathbb{H}^n)$, resp. $\textbf{w}^\varepsilon\in C^\infty(\mathbb{H}^n,H\mathbb{H}^n)$, are defined as
\begin{equation*}
	\lambda^\varepsilon(x):=\rho_\varepsilon\ast\lambda(x)=\int_{\mathbb{H}^n}\rho_\varepsilon(x\cdot y^{-1})d\lambda(y),
\end{equation*}
and
\begin{equation*}
	\textbf{w}^\varepsilon:=\sum_{j=1}^{2n}\textbf{w}^\varepsilon_jX_j,
\end{equation*}
where $\textbf{w}_j^\varepsilon(x):=\rho_\varepsilon\ast\textbf{w}_j(x)$, for any $j=1,\ldots,2n$.

\subsubsection{Riemannian approximation of $\mathbb{H}^n$}\label{riemannianapproxsubsection}
The Heisenberg group $\left(\mathbb{H}^n,d_{SR}\right)$ arises as the pointed Hausdorff-Gromov limit of Riemannian manifolds, defined as follows, in which the non-horizontal direction is increasingly penalized. One can equip $\mathbb{H}^n\equiv\mathbb{R}^{2n+1}$ with a left-invariant Riemannian metric $g_\epsilon$, with $\epsilon>0$, that has
\begin{equation*}
	X_1,\ldots,X_{2n},\epsilon X_{2n+1}
\end{equation*}
as orthonormal frame for $T\mathbb{H}^n$. We denote by $|\cdot|_\epsilon:=\sqrt{g_\epsilon(\cdot,\cdot)}$ and by $d_\epsilon$ the associated Riemannian distance: it turns out to be left-invariant with respect to \eqref{grouplaw}, since $X_1,\ldots, X_{2n}, \epsilon X_{2n+1}$ are themselves left invariant. It follows from the definition that, for fixed $x,y\in \mathbb{H}^n$, the function $d_\epsilon(x,y)$ is decreasing in $\epsilon$ and for any $\epsilon\in(0,\overline\epsilon)$
\begin{equation*}
	d_\epsilon(x,y)\leq d_{SR}(x,y).
\end{equation*}
Moreover, for any $x,y\in \mathbb{H}^n$ it holds that
\begin{equation*}
	d_{SR}(x,y)=\lim_{\epsilon\to0}d_\epsilon(x,y),
\end{equation*}
see \cite{Gromov}, and the convergence is uniform on compact subsets, see \cite{capogna2016regularity}. In \cite[Theorem 1.1]{Ge1993} the author proved that $(\mathbb{H}^n,d_\epsilon)$ converges to $(\mathbb{H}^n,d_{SR})$ in the pointed Gromov–Hausdorff sense, as $\epsilon$ goes to $0$.

In the end, if $\varphi\in C^{\infty}(\mathbb{H}^n)$ we denote by
\begin{equation*}
    \nabla_\epsilon\varphi:=\sum_{i=1}^{2n}X_i\varphi X_i+\epsilon^2X_{2n+1}\varphi X_{2n+1}
\end{equation*}

\subsection{Optimal transport theory in $\mathbb{H}^n$}
This section contains some well-known results about the Monge-Kantorovich problem in $\mathbb{H}^n$ with cost function equal to the Carnot-Carathéodory distance. 

Let $(M_1,d_1)$ and $(M_2,d_2)$ be two Polish spaces and $F:M_1\to M_2$ be a Borel map. Given a finite Radon measure $\lambda\in\mathcal{M}(M_1)$, we denote by $F_{\#}\lambda\in\mathcal{M}(M_2)$ the push-forward measure, that is the finite Radon measure defined as $F_{\#}\lambda(A):=\lambda(F^{-1}(A)), \forall A$ Borel set in $M_2$. By $\mathcal{P}(M)$ we denote the subset of $\mathcal{M}(M)$ consisting of probability measures on the metric space $M$, and by $\mathcal{P}_c(M)$ the subset of compactly supported ones. 

Now, given  $\mu,\nu\in\mathcal{P}_c(\mathbb{H}^n)$, we denote by
\begin{equation*}
	\Pi(\mu,\nu)=\big{\{}\gamma\in\mathcal{P}(\mathbb{H}^n\times\mathbb{H}^n): (\pi_1)_{\#}\gamma=\mu, (\pi_2)_{\#}\gamma=\nu\big{\}}
\end{equation*}
the set of transport plans between $\mu$ and $\nu$, where $\pi_1$ and $\pi_2$ are the projection on the first and second factor, respectively; this set is compact w.r.t. the weak convergence of measures.

It is well-known that the Monge-Kantorovich problem between $\mu$ and $\nu$, associated with the Carnot-Carathéodory distance,
\begin{equation}\label{MKH}
	\inf_{\gamma\in \Pi(\mu,\nu)} \int_{\mathbb{H}^n\times\mathbb{H}^n} d_{SR}(x,y)\,d\gamma(x,y),
\end{equation}
admits at least one solution. See for instance \cite[Theorem 1.7]{Santambrogiolibro}. We denote by 
\begin{equation*}
	\Pi_1(\mu,\nu):=\{\gamma\in\Pi(\mu,\nu):\gamma \text{ solves } \eqref{MKH}\}
\end{equation*}
the set of \textit{optimal transport plans}. It is a closed subset of the compact set $\Pi(\mu,\nu)$, w.r.t. the weak convergence of measures.

Let us denote by
\begin{equation}\label{lipset}
	\text{Lip}_1(\mathbb{H}^n, d_{SR}):=\left\{u:\mathbb{H}^n\to\mathbb{R}:|u(x)-u(y)|\leq d_{SR}(x,y),\ \forall x,y\in\mathbb{H}^n\right\}.
\end{equation}
It holds the following important theorem, see for instance \cite[Proposition 3.1]{Santambrogiolibro}.
\begin{teo}[Kantorovich duality theorem] \label{1lip_potential}
There exists $\bar{u}\in\textnormal{Lip}_1(\mathbb{H}^n,d_{SR})$ that solves 
\begin{equation}\label{11maggio}
    \sup\left\{\int_{\mathbb{H}^n} ud\left(\mu-\nu\right):u\in \textnormal{Lip}_1(\mathbb{H}^n,d_{SR})\right\}.
\end{equation}
Moreover, it holds that 
\begin{equation*}
    \min_{\gamma\in\Pi(\mu,\nu)}\int_{\mathbb{H}^n\times\mathbb{H}^n} d_{SR}(x,y)\,d\gamma(x,y) 
		=\int_{\mathbb{H}^n}\bar{u}d(\mu-\nu),
\end{equation*}
and $\gamma\in \Pi(\mu,\nu)$ is optimal if, and only if, 
\begin{equation*}
	\bar{u}(x) - \bar{u}(y) = d_{SR}(x,y)  \qquad \gamma-\text{a.e. in } \mathbb{H}^n\times\mathbb{H}^n.
\end{equation*}
\end{teo}

Any $u\in\text{Lip}_1(\mathbb{H}^n,d_{SR})$ that solves \eqref{11maggio} is a \textit{Kantorovich potential}.

If at least one between $\mu$ and $\nu$ is absolutely continuous with respect to the Haar measure, then the explicit representation of geodesics implies that any optimal transport plan is concentrated on the set $E\subset\mathbb{H}^n\times\mathbb{H}^n$ of pairs of points connected by a unique geodesic, see \eqref{KAPPA} for its definition.

\begin{lemma}\label{19marzo}
	Let us suppose that either $\mu\ll\mathcal{L}^{2n+1}$, or $\nu\ll\mathcal{L}^{2n+1}$, and $\gamma\in\Pi_1(\mu,\nu)$. Then for $\gamma$-a.e. $(x,y)$ in $\mathbb{H}^n\times\mathbb{H}^n$ there exists a unique geodesic between $x$ and $y$. That is,
 $$\gamma(\{(x,y)\in\mathbb{H}^n\times\mathbb{H}^n: x^{-1}\cdot y\in L\})=0.$$
\end{lemma}

For the proof of this result we refer to \cite[Lemma 4.1]{DePascale2}.

\section{Properties of transport rays and Kantorovich potentials in $\mathbb{H}^n$}\label{Transportrays}
In the spirit of \cite[Chapter 3]{Santambrogiolibro} and \cite[Section 2]{Feldman} we extend to the Heisenberg setting some properties of transport rays and of transport set.

\subsection{Transport rays and transport set}

Let us consider $\mu,\nu\in\mathcal{P}_c(\mathbb{H}^n)$ and a Kantorovich potential $u\in\text{Lip}_1(\mathbb{H}^n,d_{SR})$.  

Let us just remark that if $x\in\textnormal{supp}(\mu)$ and $y\in\textnormal{supp}(\nu)$ are such that 
\begin{equation}\label{inutile}
    u(x)-u(y)=d_{SR}(x,y),
\end{equation}
then 
\begin{equation*}
	u(\sigma(t))=u(x)-d_{SR}\left(x,\sigma(t)\right),\quad\forall t\in[0,1],
\end{equation*}
where $\sigma$ is any geodesic between $x$ and $y$.

\begin{deff}
A transport ray is a non-trivial geodesic $\sigma:[0,1]\to\mathbb{H}^n$ such that
\begin{enumerate}
    \item $\sigma(0)\in\textnormal{supp}(\mu)$ and $\sigma(1)\in\textnormal{supp}(\nu)$;
    \item $u(\sigma(0))-u(\sigma(1))=d_{SR}(\sigma(0),\sigma(1))$;
\end{enumerate}
\end{deff}
We denote by $\mathcal{T}_1$ the set of all points which lie on transport rays
\begin{equation*}
	\mathcal{T}_1:=\left\{x\in \sigma([0,1]):\sigma\text{ is transport ray}\right\},
\end{equation*}
and by $\mathcal{T}_0$ the complementary set of \textit{rays of length zero}
\begin{equation*}
	\mathcal{T}_0:=\Big{\{}x\in\textnormal{supp}(\mu)\cap\textnormal{supp}(\nu):|u(x)-u(z)|<d_{SR}(x,z), \forall z\in\textnormal{supp}(\mu)\cup\textnormal{supp}(\nu),x\not=z\Big{\}}.
\end{equation*}
We will call \textit{transport set} the set
\begin{equation}\label{tansportset}
	\mathcal{T}:=\mathcal{T}_1\cup\mathcal{T}_0.
\end{equation}
We observe that
\begin{equation}\label{nonemptytransportset}
	\textnormal{supp}(\mu)\cup\textnormal{supp}(\nu)\subseteq\mathcal{T}.  
\end{equation}

Moreover, as in \cite[Lemma 8]{Feldman} the following result holds.
\begin{teo}\label{comptranspset}
	The transport set $\mathcal{T}$ is compact.
\end{teo}
\begin{proof}
	Thanks to Hopf-Rinow theorem it is enough to prove that $\mathcal{T}$ is a closed and bounded set. Let us consider the function $v:\mathbb{H}^n\times\mathbb{H}^n\to\mathbb{R}$, $v(x,y)=u(x)-u(y)$. Since $v$ is continuous, it attains a maximum $L<\infty$ on $\textnormal{supp}(\mu)\times\textnormal{supp}(\nu)$, which is a compact set. Let us prove that $L\geq0$. If $\textnormal{supp}(\mu)\cap\textnormal{supp}(\nu)\not=\emptyset$, then $\forall x\in\textnormal{supp}(\mu)\cap\textnormal{supp}(\nu)$ we have that $(x,x)\in\textnormal{supp}(\mu)\times\textnormal{supp}(\nu)$ and $v(x,x)=0$. Otherwise, from \eqref{nonemptytransportset} it follows that $\mathcal{T}_1\not=\emptyset$, hence there exists at least one transport ray $\sigma$. If $x=\sigma(0)$ and $y=\sigma(1)$, then $v(x,y)=d_{SR}(x,y)>0$. Hence $L\geq0$.
	
	We can suppose that $A:=\mathcal{T}\setminus \textnormal{supp}(\mu)\cup\textnormal{supp}(\nu)\not=\emptyset$; otherwise, from the previous theorem it follows that $\mathcal{T}=\textnormal{supp}(\mu)\cup\textnormal{supp}(\nu)$, which is compact. Hence, any $z\in A$ lies on a transport ray $\sigma_z$. Let us denote by $a=\sigma_z(0)$ and $b=\sigma_z(1)$, then
	\begin{equation*}
		d_{SR}(a,z)+d_{SR}(b,z)=d_{SR}(a,b)=v(a,b)\leq L.
	\end{equation*}
	Hence, $A$ lies in the union of the $L$-neighborhoods of the compact sets $\textnormal{supp}(\mu)$ and $\textnormal{supp}(\nu)$, thus $\mathcal{T}$ is bounded. 
	
	Let us prove that $\mathcal{T}$ is closed. Let us consider $(z_n)_{n\in\mathbb{N}}\subseteq\mathcal{T}$, converging to some $z$, we prove that $z\in\mathcal{T}$. If there exists a subsequence $(z_{n_k})_{k\in\mathbb{N}}\subseteq \textnormal{supp}(\mu)\cup\textnormal{supp}(\nu)$, then $z\in \textnormal{supp}(\mu)\cup\textnormal{supp}(\nu)$ by compactness of $\textnormal{supp}(\mu)$ and $\textnormal{supp}(\nu)$. Let us suppose that $z_n\in A, \forall n\in\mathbb{N}$: there exists a transport ray $\sigma_n$, whose endpoints we denote by $a_n=\sigma_n(0)$ and $b_n:=\sigma_n(1)$. There exist two subsequences 
	\begin{equation}\label{transportset1}     a_{n_j}\to a\in\textnormal{supp}(\mu)\text{ and } b_{n_j}\to b\in\textnormal{supp}(\nu), \text{ when }  j\to\infty.
	\end{equation}
	Since 
	\begin{equation*}
		d_{SR}(z_{n_j},a_{n_j})+d_{SR}(z_{n_j},b_{n_j})=d_{SR}(a_{n_j},b_{n_j})=u(a_{n_j})-u(b_{n_j}),\quad \forall n_j,
	\end{equation*}
	we have that 
	\begin{equation}\label{transportset2}
		d_{SR}(z,a)+d_{SR}(z,b)=d_{SR}(a,b)=u(a)-u(b).
	\end{equation}
	Hence there are two possibilities
	\begin{enumerate}
		\item $a=b\underset{\eqref{transportset2}}{\Longrightarrow}z=a=b\in\mathcal{T};$
		\item $a\not=b\underset{\eqref{transportset1}+\eqref{transportset2}}{\Longrightarrow}z\in\mathcal{T}_1.$
	\end{enumerate}
\end{proof}

\subsection{Differentiability of Kantorovich potential}
Following \cite[Lemma 10]{Feldman}, one can prove that any Kantorovich potential $u$ is Pansu differentiable in the interior of transport rays.
\begin{prop}\label{diffu} Let $u\in\textnormal{Lip}_1(\mathbb{H}^n, d_{SR})$ and $x,y\in\mathbb{H}^n, x\not=y$ such that $u(x)-u(y)=d_{SR}(x,y)$. Let $\sigma:[0,1]\to\mathbb{H}^n$ be a geodesic between $x$ and $y$, starting from $x$. Then $u$ is Pansu differentiable at $\sigma(t)$, for all $t\in]0,1[$, and
	\begin{equation}
		\nabla_Hu(\sigma(t))=-\frac{\dot\sigma(t)}{|\dot\sigma(t)|_H}.
	\end{equation}
\end{prop}

For the proof of this result we need the following two lemmas. The first one collects a differentiability property of the Carnot-Charatheodory distance function from a fixed point $y\in \mathbb{H}^n$. Let us denote $L_y:= y\cdot L$.

\begin{lemma} \label{prop-distcc} 

Let $y\in\mathbb{H}^n$, then function $d_y(\cdot) := d_{SR}(\cdot,y)$ is of class $C^{\infty}$ in the euclidean sense on $\mathbb{H}^n \backslash L_y$. In particular $d_y$ is Pansu differentiable on $\mathbb{H}^n \backslash L_y$. 

Moreover, if $x\in\mathbb{H}^n\setminus L_y$ and $\sigma:[0,1]\to\mathbb{H}^n$ is the geodesic between $x$ and $y$, starting from $y$, then
\begin{equation*}
	\nabla_Hd_y(x)=\frac{\dot{\sigma}(1)}{|\dot{\sigma}(1)|_H}. 
\end{equation*}
\end{lemma}

\begin{proof} 
The core of the proof is contained in \cite[Lemma 3.11]{Rigot}.

If we denote by $\Phi(\chi,\varphi):=y\cdot\sigma_{\chi,\theta}(1)$, where $\sigma_{\chi,\theta}:[0,1]\to\mathbb{H}^n$ is any geodesic starting from $0$ as in Theorem~\ref{geod}, then this map is a $C^\infty$-diffeomorphism from $\mathbb{R}^{2n}\setminus\{0\} \times (-2\pi,2\pi)$ onto $\mathbb{H}^n\setminus L_y$ in the usual euclidean sense. If $z=\Phi(\chi,\theta)\in\mathbb{H}^n\setminus L_y$, then $d_y(z)=|\chi|$. Therefore $d_y$ is $C^\infty$ in $\mathbb{H}^n\setminus L_y$, thus it Pansu differentiable in $\mathbb{H}^n\setminus L_y$. In addition, the condition $d_y\in \text{Lip}_1(\mathbb{H}^n,d_{SR})$ implies that $|\nabla_Hd_y(z)|\leq1, \forall z\in \mathbb{H}^n\setminus L_y$. 

Now let $x\in\mathbb{H}^n\setminus L_y$ and $\sigma:[0,1]\to\mathbb{H}^n$ be the geodesic between $x$ and $y$. Since $d_y(\sigma(t))=td_{SR}(x,y)$ for all $t\in[0,1]$, and $\sigma(t)\in \mathbb{H}^n\setminus L_y$ for all $t\in]0,1]$, we can differentiate w.r.t. $t$:
\begin{multline}\label{nonso}
	d_{SR}(x,y)=\frac{d}{dt}d_y(\sigma(t))=\sum_{j=1}^{2n}X_j(d_y(\sigma(t)))\dot\sigma_j(t)\leq |\nabla_Hd_y(\sigma(t))|_H|\dot\sigma(t)|_H\leq d_{SR}(x,y),
\end{multline}
where we used the fact that $\sigma$ is a geodesic, then $|\dot\sigma(t)|_H=d_{SR}(x,y)$. Hence, all the inequalities in \eqref{nonso} are equalities: in particular we get that 
\begin{equation*}
    \nabla_Hd_y(\sigma(t))=\frac{\dot{\sigma}(t)}{|\dot{\sigma}(t)|_H},\quad \forall t\in]0,1].
\end{equation*}
\end{proof}

\begin{lemma}\label{Pansudiffprop}
	Lef $f,g,h:\mathbb{H}^n\to\mathbb{R}$ three functions such that $f(x)\leq g(x)\leq h(x)$, for all $x\in\mathbb{H}^n$. Let $y\in\mathbb{H}^n$ such that $f(y)=g(y)=h(y)$ and $f,h$ are Pansu differentiable at $y$, with $\nabla_H f(y)=\nabla_Hh(y)$. Then, $g$ is Pansu differentiable in $y$ and $\nabla_Hg(y)=\nabla_H f(y)=\nabla_Hh(y)$.
\end{lemma}

\begin{proof}
Let $y\in\mathbb{H}^n$ as in the hypothesis and $x$ be any point in $\mathbb{H}^n$. From \eqref{9maggio} we have that
\begin{equation*}
	D_Hf(y)(x)=\left\langle\nabla_Hf(y),\pi_{y}(x)\right\rangle_H=\left\langle\nabla_Hh(y),\pi_{y}(x)\right\rangle_H =D_Hh(y)(x).
\end{equation*}
Therefore,
\begin{equation*}
\aligned
	f(x)-f(y)-D_Hf(y)(y^{-1}\cdot x)&\leq g(x)-g(y)-\underbrace{D_Hf(y)(y^{-1}\cdot x)}_{=D_Hh(y)(y^{-1}\cdot x)}\\&\leq h(x)-h(y)-D_Hh(y)(y^{-1}\cdot x).
\endaligned
\end{equation*}
If we divide the previous inequalities by $d_{SR}(x,y)$ and we let $x$ tend to $y$, by using Pansu differentiability of $f$ and $g$ and \cite[Proposition 5.6]{Franchi1} again, we get the thesis.
\end{proof}

\begin{proof}[Proof of Proposition \ref{diffu}]
Let $t_0\in]0,1[$ and $a,b\in\sigma([0,1])$ such that $u(a)>u(\sigma(t_0))>u(b)$. 

The description of geodesics given in Theorem \ref{geod} implies that $\sigma(s)^{-1}\cdot\sigma(t)\not\in L$, for any $s,t\in]0,1[$. Therefore, it follows from Lemma \ref{prop-distcc} that the functions $d_{a}(\cdot)$ and $d_{b}(\cdot)$ are smooth in a neighborhood of $\sigma(t_0)$. Since $u\in\textnormal{Lip}_1(\mathbb{H}^n,d_{SR})$, it holds that $$u(a)-u(z)\leq d_{a}(z),\quad \forall z\in\mathbb{H}^n.$$ Moreover $a$ and $b$ lie on the geodesic between $x$ and $y$, then $u(a)=u(b)+d_{SR}(a,b)$, and hence
\begin{equation*}
	d_{b}(z)\geq u(z)-u(b)\geq d_{SR}(a,b)-d_{a}(z),\quad \forall z\in\mathbb{H}^n,
\end{equation*}
where equalities hold if $z=\sigma(t_0)$, since $\sigma$ is a geodesic. Moreover, from Lemma \ref{prop-distcc} again, it follows that $\nabla_Hd_{b}(\sigma(t_0))=-\nabla_H d_{a}(\sigma(t_0))=-\frac{\dot\sigma(t_0)}{|\dot\sigma(t_0)|_H}$. Hence, from Lemma \ref{Pansudiffprop}, it follows that $u$ is Pansu differentiable in $\sigma(t_0)$ and
\begin{equation*}
    \nabla_Hu(\sigma(t_0))=-\frac{\dot\sigma(t_0)}{|\dot\sigma(t_0)|_H}.
\end{equation*}
\end{proof}

\section{The Beckmann problem for $p=1$}\label{Beckmannp1}

\subsection{Vector transport density in $\mathbb{H}^n$}
In the spirit of \cite[Chapter3]{Santambrogiolibro} we define a vector version of the horizontal transport density, introduced \cite{circelli2023transport}.

From now on, we fix a selection of geodesics 
\begin{equation}\label{18marzo2025}
\aligned
	S:\mathbb{H}^n\times\mathbb{H}^n&\to \text{Geo}(\mathbb{H}^n)\\
  (x,y) &\mapsto S(x,y)=\sigma_{x,y}
    \endaligned
\end{equation}
that is $\gamma$-measurable for any $\gamma\in\Pi(\mu,\nu)$, see \eqref{19marzo1}. Note that throughout the rest of this section, all results depend a priori on the measurable selection $S$ of geodesics. However, this is not a serious issue for our purposes.

Hence, one may associate with any $\gamma\in\Pi_1(\mu,\nu)$ a positive and finite Radon measure $a_\gamma\in\mathcal{M}_+(\mathbb{H}^n)$, called horizontal transport density, defined as
\begin{equation*}
	\int_{\mathbb{H}^n}\varphi\, da_\gamma:=\int_{\mathbb{H}^n\times\mathbb{H}^n}\left(\int_0^1\varphi(\sigma_{x,y}(t))|\dot{\sigma}_{x,y}(t)|_Hdt\right)d\gamma(x,y),
\end{equation*}
for any $\varphi\in C(\mathbb{H}^n)$.
If either $\mu\ll\mathcal{L}^{2n+1}$, or $\mu\ll\mathcal{L}^{2n+1}$, then $\gamma(\mathbb{H}^n\times\mathbb{H}^n\setminus E)=0$ and therefore the definition of $a_\gamma$ is independent of the choice of $S$. 

One may also associate with any $\gamma\in \Pi_1(\mu,\nu)$ a vector version of  the transport density $a_\gamma$.

\begin{deff}[Vector horizontal transport density]\label{defw}
	For any $\gamma\in\Pi_1(\mu,\nu)$, the \textit{vector horizontal transport density associated with $\gamma$} is the $H\mathbb{H}^n$-valued finite Radon measure $\textnormal{\textbf{w}}_\gamma\in\mathcal{M}(\mathbb{H}^n,H\mathbb{H}^n)$ defined as
	\begin{equation}\label{horizontaltransportflow}
		\int_{\mathbb{H}^n}\phi \cdot d\textnormal{\textbf{w}}_{\gamma}:=\int_{\mathbb{H}^n\times\mathbb{H}^n}\left(\int_0^1\langle\phi\left(\sigma_{x,y}(t)\right),\dot{\sigma}_{x,y}(t)\rangle_Hdt\right)d\gamma(x,y),
	\end{equation}
	for any $\phi\in C(\mathbb{H}^n,H\mathbb{H}^n)$.
\end{deff}

We defined the scalar transport density, resp. the vector one, in duality with $C(\mathbb{H}^n)$, resp. $C(\mathbb{H}^n,H\mathbb{H}^n)$ because they are both supported on the transport set $\mathcal{T}$, which is compact, see Theorem \ref{comptranspset}. 

In particular, if $\overline\gamma\in\Pi_1(\mu,\nu)$, we have
\begin{equation}\label{densityflow}
	|\textnormal{\textbf{w}}_{\overline\gamma}|\leq a_{\overline\gamma},
\end{equation}
as measures. Therefore, it follows from the very definition of $a_{\overline\gamma}$, and the fact that $\overline\gamma$ is an optimal transport plan, that
\begin{equation}\label{totalvariationtransportflow}
	\|\textbf{w}_{\overline\gamma}\|_{\mathcal{M}(\mathbb{H}^n,H\mathbb{H}^n)}\leq \min_{\gamma\in\Pi(\mu,\nu)}\int_{\mathbb{H}^n\times\mathbb{H}^n}d_{SR}(x,y)d\gamma(x,y).
\end{equation}

\begin{prop}
    If $u\in\textnormal{Lip}(\mathbb{H}^n,d_{SR})$ is a Kantorovich potential, then for any $\gamma\in\Pi_1(\mu,\nu)$, the measure $|\textnormal{\textbf{w}}_\gamma|$ is absolutely continuous w.r.t. $a_\gamma$, with density $-\nabla_Hu$  $$\textnormal{\textbf{w}}_{\gamma}=-(\nabla_Hu)a_\gamma.$$
\end{prop}
\begin{proof}
    Let $\gamma\in\Pi_1(\mu,\nu)$. Then, Proposition \ref{diffu} implies that 
\begin{equation*}
	\dot{\sigma}_{x,y}(t)=|\dot{\sigma}_{x,y}(t)|_H\frac{\dot{\sigma}_{x,y}(t)}{|\dot{\sigma}_{x,y}(t)|_H}=-d_{SR}(x,y)\nabla_Hu(\sigma_{x,y}(t)),
\end{equation*}
for every $t\in]0,1[$ and for $\gamma$-almost every $(x,y)$ (with $x\neq y$ otherwise both expressions vanish), where  $\sigma_{x,y}=S(x,y)$ is the geodesic between $x$ and $y$, selected by the map $S$ in \eqref{18marzo2025}.
Hence, it follows that 
\begin{align*}
	\int_{\mathbb{H}^n}\phi\cdot d\textnormal{\textbf{w}}_{\gamma}&=\int_{\mathbb{H}^n\times\mathbb{H}^n}\left(\int_0^1-d_{SR}(x,y)\left\langle\nabla_H u(\sigma_{x,y}(t)),\phi(\sigma_{x,y}(t)) \right\rangle_Hdt\right)d\gamma(x,y)=\\&=-\int_0^1\left(\int_{\mathbb{H}^n\times\mathbb{H}^n} \left\langle\nabla_H u(\sigma_{x,y}(t)),\phi(\sigma_{x,y}(t)) \right\rangle_Hd_{SR}(x,y)d\gamma(x,y) \right)dt=\\
	&=-\int_0^1\left(\int_{\mathbb{H}^n} \left\langle\nabla_H u(z),\phi(z) \right\rangle_Hd(S_t)_\#(d_{SR}\gamma)(z) \right)dt,
\end{align*}
for every $\phi\in C(\mathbb{H}^n,H\mathbb{H}^n)$, where $S_t=e_t\circ S$. Analogously
\begin{equation*}
	\int_{\mathbb{H}^n}\varphi da_\gamma=\int_0^1\left(\int_{\mathbb{H}^n}\varphi(z) d(S_t)_\#(d_{SR}\gamma)(z)\right)dt,
\end{equation*}
for every $\varphi\in C(\mathbb{H}^n)$. Since, by definition,  $\textbf{w}_\gamma$ and $a_\gamma$ are concentrated on the set of differentiability points of the Kantorovich potential $u$,  it follows that
\begin{equation*}
	\int_{\mathbb{H}^n}\phi\cdot d\textnormal{\textbf{w}}_{\gamma}=-\int_{\mathbb{H}^n}\left\langle\nabla_Hu,\phi\right\rangle_Hda_\gamma=-\int_{\mathbb{H}^n}\phi\cdot d\left(\nabla_Hu\right)a_\gamma,
\end{equation*}
for every $\phi\in C(\mathbb{H}^n,H\mathbb{H}^n)$; hence $-\nabla_Hu$ is the density of the measure $|\textbf{w}_{\gamma}|$ w.r.t. $a_\gamma$, and we write
\begin{equation*}
	\textnormal{\textbf{w}}_{\gamma}=-(\nabla_Hu)a_\gamma.
\end{equation*}
\end{proof}

From \eqref{densityflow} and \cite[Theorem 3.3 and Theorem 3.6]{circelli2023transport} the next theorem easily follows. 
\begin{teo}\label{4giugno}[Absolute continuity and summability of vector horizontal transport densities]
	If $\mu\ll\mathcal{L}^{2n+1}$, either $\nu\ll\mathcal{L}^{2n+1}$, then there exists $\gamma\in\Pi_1(\mu,\nu)$ such that $\textnormal{\textbf{w}}_{\gamma}\in L^1(\mathbb{H}^n,H\mathbb{H}^{n})$.
    Moreover, if $\mu\in L^p$ for some $p\in[1,\infty]$, it holds: if $p<\frac{2n+3}{2n+2}$ then there exists $\gamma\in\Pi_1(\mu,\nu)$ such that $\textnormal{\textbf{w}}_{\gamma}\in L^p(\mathbb{H}^n,H\mathbb{H}^n)$; otherwise there exists $\gamma\in\Pi_1(\mu,\nu)$ such that $\textnormal{\textbf{w}}_{\gamma}\in L^s(\mathbb{H}^n,H\mathbb{H}^n)$ for $s<\frac{2n+3}{2n+2}$.
\end{teo}

\subsection{The problem}
Given any optimal transport plan $\gamma\in\Pi_1(\mu,\nu)$ we can test the vector horizontal transport density $\textbf{w}_{\gamma}$ against $\phi=\nabla_H\varphi$, for any $\varphi\in C^\infty(\mathbb{H}^n)$. From \eqref{14marzo} it follows that 
\begin{multline}\label{divergenceoftransportflow}
	\int_{\mathbb{H}^n}\nabla_H\varphi\cdot d{\textbf{w}}_{\gamma}=\int_{\mathbb{H}^n\times\mathbb{H}^n}\left(\int_0^1\langle\nabla_H\varphi\left(\sigma_{x,y}(t)\right),\dot\sigma_{x,y}(t) \rangle_Hdt\right)d\gamma(x,y)=\\=\int_{\mathbb{H}^n\times\mathbb{H}^n}\left(\int_0^1\frac{d}{dt}\left[\varphi(\sigma_{x,y}(t))\right]dt\right)d\gamma(x,y)
	\\=\int_{\mathbb{H}^n\times\mathbb{H}^n}(\varphi(y)-\varphi(x))d\gamma(x,y)=-\int_{\mathbb{H}^n}\varphi d(\mu-\nu).
\end{multline}

Now, given a $H\mathbb{H}^n$-valued compactly supported Radon measure $\textbf{w}\in\mathcal{M}_c(\mathbb{H}^n,H\mathbb{H}^n)$ we can define its \textit{horizontal divergence} $\textnormal{div}_H\textbf{w}$ in duality with the smooth functions
$$\langle \textnormal{div}_H\textbf{w},\varphi\rangle:=-\int_{\mathbb{H}^n} \nabla_H\varphi(x) \cdot d\textbf{w}, \quad\forall \varphi\in C^\infty(\mathbb{H}^n).$$

With this definition in mind we can rewrite \eqref{divergenceoftransportflow} in the following way:
\begin{equation*}
	\textnormal{div}_H\textbf{w}_{\gamma}=\mu-\nu,\quad \forall \gamma\in\Pi_1(\mu,\nu),
\end{equation*}
i.e. the horizontal divergence of the measure $\textbf{w}_{\gamma}$ is the signed Radon measure $\mu-\nu$, $\forall \gamma\in\Pi_1(\mu,\nu)$.
This says that the following infimum is finite,
\begin{equation}\label{beckmann}\tag{BP}
	\inf\left\{\|\textnormal{\textbf{w}}\|_{\mathcal{M}(\mathbb{H}^n,H\mathbb{H}^n)}:\textnormal{\textbf{w}}\in\mathcal{M}_c(\mathbb{H}^n,H\mathbb{H}^n),\, \textnormal{div}_H\textnormal{\textbf{w}}=\mu-\nu\right\}.
\end{equation}
This is the Heisenberg version of the well-known \textit{Beckmann problem} and it will turn out to be the limit case of the wider class of problems introduced in Section 5, as explained in the introduction of the paper.

In the spirit of \cite[Subsection 4.2]{Santambrogiolibro}, we show that this problem is equivalent to the Monge-Kantorovich problem \eqref{MKH}. 

\begin{teo}\label{7dicembre1}
The problem \eqref{beckmann} admits a solution. Moreover,
\begin{equation}
	\eqref{beckmann}=\eqref{MKH}
\end{equation}
where 
\begin{equation*}
	\eqref{MKH}=\min_{\gamma\in\Pi(\mu,\nu)}\int_{\mathbb{H}^n\times\mathbb{H}^n}d_{SR}(x,y)d\gamma(x,y),
\end{equation*}
and a solution to \eqref{beckmann} can be built from a solution to $\eqref{MKH}$, as in \eqref{horizontaltransportflow}. 
\end{teo}
\begin{proof}
First we prove the equality between \eqref{MKH} and \eqref{beckmann}. 
We start by proving that \eqref{beckmann}$\geq$\eqref{MKH}. Take an arbitrary function $\varphi\in C^\infty\cap\textnormal{Lip}_1(\mathbb{H}^n,d_{SR})$. Pansu's Theorem implies that $\|\nabla_H\varphi\|_\infty\leq 1$, see \cite{Pansu}, hence for any $\textnormal{\textbf{w}}$ admissible we get by definition that
\begin{equation*}
	\|\textnormal{\textbf{w}}\|_{\mathcal{M}(\mathbb{H}^n,H\mathbb{H}^n)}\geq\int_{\mathbb{H}^n}(-\nabla_H\varphi)\cdot d\textnormal{\textbf{w}}=\int_{\mathbb{H}^n}\varphi d(\mu-\nu).
\end{equation*}
Now we take $\varphi^\varepsilon=\rho_\varepsilon\ast u$, where $u$ is a Kantorovich potential and $\rho_\varepsilon$ is a mollifier for the group structure. It follows that $\varphi^\varepsilon\in C^\infty\cap\textnormal{Lip}_1(\mathbb{H}^n,d_{SR})$ and it converges uniformly on compact sets to the Kantorovich potential $u$, see \cite[Proposition 2.14]{comi2020gauss}. By the previous inequality and letting $\varepsilon\to0$ we get that 
\begin{equation*}
	\|\textnormal{\textbf{w}}\|_{\mathcal{M}(\mathbb{H}^n,H\mathbb{H}^n)}\geq\int_{\mathbb{H}^n}u d(\mu-\nu)=\eqref{MKH},
\end{equation*}
where we used Theorem \ref{1lip_potential} in the last equality. Since the previous inequality holds for any admissible $\textnormal{\textbf{w}}$, we may take the minimum in the left hand-side and get $$\eqref{beckmann}\geq\eqref{MKH}.$$
	
Now we will prove the converse inequality: for any $\gamma\in\Pi_1(\mu,\nu)$, we know that the vector Radon measure $\textnormal{\textbf{w}}_{\gamma}$ in $H\mathbb{H}^n$, defined in \eqref{horizontaltransportflow}, is a compactly supported measure, which satisfies the divergence constraint thanks to \eqref{divergenceoftransportflow}. Moreover from \eqref{totalvariationtransportflow} we know that
\begin{equation*}
	\|\textbf{w}_{\gamma}\|_{\mathcal{M}(\mathbb{H}^n,H\mathbb{H}^n)}\leq\eqref{MKH}.
\end{equation*}
Hence
\begin{equation*}
	\eqref{beckmann}\leq\|\textbf{w}_{\gamma}\|_{\mathcal{M}(\mathbb{H}^n,H\mathbb{H}^n)}\leq\eqref{MKH}.
\end{equation*}
\end{proof}

\begin{cor}
If either $\mu\ll\mathcal{L}^{2n+1}$, or $\nu\ll\mathcal{L}^{2n+1}$, then there exists $\textnormal{\textbf{w}}\in L^1(\mathbb{H}^n,H\mathbb{H}^{n})$ that solves \eqref{beckmann}. 
\end{cor}
\begin{proof}
The thesis follows from Theorem \ref{4giugno} and Theorem \ref{7dicembre1}.
\end{proof}

As a consequence of Theorem \ref{7dicembre1}, one also obtains the following dual formulation.

\begin{cor}
It holds that 
\begin{equation}\label{dualmongechapter2}\tag{DP}
	\max\left\{\int_{\mathbb{H}^n}ud(\mu-\nu):u\in HW^{1,\infty},\|\nabla_Hu\|_\infty\leq1\right\}=\eqref{beckmann}.
\end{equation}
\end{cor}

\begin{proof}
The thesis follows from Theorem \ref{1lip_potential} and Theorem \ref{7dicembre1}.
\end{proof}

For any $\gamma\in\Pi_1(\mu,\nu)$ and any Kantorovich potential $u\in\textnormal{Lip}_1(\mathbb{H}^n,d_{SR})$, the pair $(a_\gamma,u)$ solves the Monge-Kantorovich system
\begin{equation*}
	\begin{cases}
		\textnormal{div}_H((\nabla_Hu)a_\gamma)=\mu-\nu,\\
		|\nabla_Hu|_H=1,\quad a_\gamma-a.e.,\\
		|\nabla_Hu|_H\leq1.
	\end{cases}
\end{equation*}
The second condition holds because $|\nabla_Hu|_H=1$ in the relative interior of transport rays and $a_\gamma$ is supported on the transport set. 

\subsection{Lagrangian formulation}
 
Following \cite[Lecture 9, section 3]{ambrosio2021lectures},  we now introduce a Lagrangian formulation of \eqref{beckmann} in $\mathbb{H}^n$,  Let consider the space $C([0,1],\mathbb{R}^{2n+1})$ equipped with the topology of uniform convergence and let us denote by $\mathcal{P}(C([0,1],\mathbb{R}^{2n+1}))$ the set of Borel probability measures over this space. We are interested in  the problem 
\begin{equation}\label{Dyn}
	\inf\Bigg{\{}\int_{C([0,1],\mathbb{R}^{2n+1})}\ell(\sigma)dQ(\sigma): Q\in\mathcal{P}(C([0,1],\mathbb{R}^{2n+1})),(e_0)_\#Q=\mu, (e_1)_\#Q=\nu\Bigg{\}},
\end{equation}
where the functional $\ell$ is the length induced by the sub-Riemannian distance, 
\begin{equation}\label{lengthstruc}
	\ell(\sigma):=\sup\left\{ \sum_{i=1}^{k}d_{SR}(\sigma(t_{i-1}),\sigma(t_i)) : k\in\mathbb{N}, 0=t_0<t_1<\ldots<t_{k-1}<t_{k}=1\right\}.
\end{equation}
This definition of length makes sense even for continuous curves and it extends the sub-Riemannian one, in the sense that if $\sigma$ is horizontal, then $\ell(\sigma)=\ell_{SR}(\sigma)$. See for instance \cite[Theorem 1.3.5]{monti2001distances}. It immediately follows that the functional $C([0,1],\mathbb{R}^{2n+1})\ni\sigma\mapsto\ell(\sigma)\in[0,+\infty]$ is lower semicontinuous w.r.t. the topology of uniform convergence, hence Borel.

Hence, the following theorem holds.
\begin{teo}\label{dynteo}
If $\mu,\nu\in\mathcal{P}_c(\mathbb{H}^n)$, then \eqref{Dyn} admits a solution. 
	
Moreover, it holds
\begin{equation*}
	\eqref{Dyn}=\min_{\gamma\in\Pi(\mu,\nu)}\left\{\int_{\mathbb{H}^n\times\mathbb{H}^n}d_{SR}(x,y)d\gamma(x,y)\right\}<+\infty
\end{equation*}
and $Q$ is optimal if and only if $Q$ is supported on the set of minimizing horizontal curves and $(e_0,e_1)_\#Q\in\Pi_1(\mu,\nu)$.
\end{teo}

\begin{proof}
Let us first observe that, for any $\gamma \in \Pi(\mu, \nu)$, the measure 
\begin{equation*}
	Q_\gamma := S_\# \gamma \in \mathcal{P}(C([0,1], \mathbb{R}^{2n+1})),
\end{equation*}
satisfies 
\begin{equation*}
    (e_0)_\# Q = \mu, \quad (e_1)_\# Q = \nu,
\end{equation*}
by definition. This ensures that the set of admissible measures $Q$ for \eqref{Dyn} is nonempty. Moreover, $Q_\gamma$ is supported on Geo($\mathbb{H}^n$) by definition.

Let $Q \in \mathcal{P}(C([0,1], \mathbb{R}^{2n+1}))$ be admissible for \eqref{Dyn}. Then,
\begin{equation}\label{9.9}
\aligned
	\int_{C([0,1], \mathbb{R}^{2n+1})} \ell(\sigma) \, dQ(\sigma) &\geq \int_{C([0,1], \mathbb{R}^{2n+1})} d_{SR}(\sigma(0), \sigma(1)) \, dQ(\sigma)\\
	&=\int_{\mathbb{H}^n \times \mathbb{H}^n} d_{SR}(x, y) \, d (e_0, e_1)_\# Q(x, y)\\ &\geq \min_{\gamma \in \Pi(\mu, \nu)} \left\{ \int_{\mathbb{H}^n \times \mathbb{H}^n} d_{SR}(x, y) \, d\gamma(x, y) \right\},
\endaligned
\end{equation}
where the first inequality follows from the definition of $\ell$ \eqref{lengthstruc}, the second equality follows from the definition of push-forward, and the last inequality follows from the fact that $(e_0, e_1)_\# Q \in \Pi(\mu, \nu)$.

To prove the converse inequality, let us consider an optimal transport plan $\overline{\gamma} \in \Pi_1(\mu, \nu)$. From the discussion at the beginning of the proof, the measure $Q_{\overline{\gamma}} := S_\# \overline{\gamma} \in \mathcal{P}(C([0,1], \mathbb{R}^{2n+1}))$ is admissible for \eqref{Dyn} and supported on Geo($\mathbb{H}^n$). Hence,
\begin{equation}\label{9.10}
\aligned
	\min_{\gamma \in \Pi(\mu, \nu)} \left\{ \int_{\mathbb{H}^n \times \mathbb{H}^n} d_{SR}(x, y) \, d\gamma(x, y) \right\} &= \int_{\mathbb{H}^n \times \mathbb{H}^n} d_{SR}(x, y) \, d\overline{\gamma}(x, y) \\
	&= \int_{\text{Geo}(\mathbb{H}^n)} d_{SR}(\sigma(0), \sigma(1)) \, dQ_{\overline{\gamma}}(\sigma)\\& = \int_{\text{Geo}(\mathbb{H}^n)} \ell(\sigma) \, dQ_{\overline{\gamma}}(\sigma).
\endaligned
\end{equation}
Thus, the equality of optimization problems follows. Moreover, if $Q \in \mathcal{P}(C([0,1], \mathbb{R}^{2n+1}))$, then $Q$ is an optimizer for \eqref{Dyn} if and only if \eqref{9.9} holds with equalities, which occurs if and only if $Q$ is supported on the set of minimizing horizontal curves and $(e_0, e_1)_\# Q \in \Pi_1(\mu, \nu)$. 
\end{proof}

This theorem is independent on the compactness of supports and on the ambient space: it holds in more general Polish metric spaces $(X,d)$ with $\mu,\nu\in\mathcal{P}(X)$ such that $$\int_{X}d(x,0)d\mu(x)+\int_{X}d(0,y)(y)d\nu(y)<+\infty,$$
and $\ell$ is the length induced by the distance $d$ as in \eqref{lengthstruc}.

\section{The Beckmann problem for $1<p<+\infty$}\label{Beckmannpgeq1}

From now on we suppose $\Omega$ is a bounded domain, with regular $C^{1,1}$ boundary (in the Euclidean sense). Such a $\Omega$ support a $q$-Poincaré inequality, i.e. there exists $c=c(n,q,\Omega)>0$ such that
\begin{equation}\label{poincaréineq}
	\int_\Omega|\varphi(x)-\varphi_\Omega|^qdx\leq c\int_\Omega|\nabla_H\varphi(x)|^qdx,\quad \forall\varphi\in HW^{1,q}(\Omega)
\end{equation} 
for any $1\leq q<+\infty$, where $\varphi_\Omega=\frac{1}{\mathcal{L}^{2n+1}(\Omega)}\int_\Omega \varphi\ dx$. See \cite[(1.26) and Theorem 1.21]{garofalo1996isoperimetric}.

We consider the space
\begin{equation*}
    H:=\left\{\sigma\in AC([0,1],\overline{\Omega}):\sigma\text{ is horizontal}\right\},
\end{equation*}
viewed as subset of $C([0,1],\overline{\Omega})$ equipped with the topology of uniform convergence. We stress the fact that the regularity assumption on the boundary of $\Omega$ guarantees that any two points in $\overline\Omega$ can be connected by an element of $H$, see \cite[Lemma 4.1]{circelli2023transport}. This assumption replaces the stronger convexity assumption typically used in the Euclidean setting, since there are no non-trivial geodesically convex set in $\mathbb{H}^n$, see \cite{monti2005geodetically}.

\subsection{Vector traffic intensities in $\mathbb{H}^n$}
In this subsection we will see how also the congested optimal transport problem is related with some minimization problems, over a set of vector fields with divergence constraint.

Let $\mu,\nu\in\mathcal{P}(\overline{\Omega})$ be two probability measures over the closure of $\Omega$ and let us consider the space $C([0,1],\overline{\Omega})$ equipped with the topology of uniform convergence. A horizontal traffic plan between $\mu$ and $\nu$ is a Borel probability measure $Q\in\mathcal{P}(C([0,1],\overline{\Omega}))$ such that 
\begin{enumerate}
    \item $\int_{C([0,1],\overline{\Omega})}\ell(\sigma)dQ(\sigma)<+\infty;$
    \item $Q$ is concentrated on the set $H$;
    \item $(e_0)_\#Q=\mu$ and $(e_1)_\#Q=\nu$.
\end{enumerate}
We denote by
\begin{equation*}
	\mathcal{Q}_H(\mu,\nu):=\left\{Q\in\mathcal{P}(C([0,1],\overline{\Omega})): Q\textnormal{ horizontal traffic plan between }\mu\textnormal{ and }\nu\right\}.
\end{equation*}
As for the transport plans, one can also associate with any horizontal traffic plan $Q\in\mathcal{Q}_H(\mu,\nu)$ both a scalar and a vector measures. The scalar measure is the horizontal traffic intensity $i_Q\in\mathcal{M}_+(\overline\Omega)$ defined as 
\begin{equation*}
	\int_{\overline\Omega} \varphi(x) d i_{Q}(x):=\int_{C([0,1], \overline\Omega)}\left(\int_0^1\varphi(\sigma(t))|\dot{\sigma}(t)|_Hdt\right) d Q(\sigma),\quad \forall \varphi \in C(\overline\Omega).
\end{equation*}

\begin{exa}\label{22marzo1}
Let us consider the two measures
\begin{equation*}
	\mu:=\delta_{x}\text{ and }\nu:=\delta_{y}, \text{ for some }x\in\mathbb{H}^n \text{ and some }y\in\mathbb{H}^n\setminus L_x,
\end{equation*}
and let $\sigma:[0,1]\to\mathbb{H}^n$ be the unique geodesic between $x$ and $y$, according to Theorem \ref{geod}. Let us consider a bounded domain $\Omega$, with regular $C^{1,1}$ boundary, such that $x,y\in\Omega$ and $\sigma([0,1])\subset\Omega$. Then, the measure
	\begin{equation*}
		Q:=\delta_{\sigma}\in\mathcal{Q}_H(\mu,\nu),
	\end{equation*}
	and its horizontal traffic intensity is precisely
	\begin{equation*}
		i_Q=\mathcal{H}^1\llcorner_{\sigma([0,1])}.
	\end{equation*}
\end{exa}

As for the vector measure, we introduce now the following definition.

\begin{deff}[Vector horizontal traffic intensity]
	Let $Q\in\mathcal{Q}_H(\mu,\nu)$ be a horizontal traffic plan admissible between $\mu$ and $\nu$. One can associate with $Q$ a $H\Omega$-valued finite Radon measure $\textnormal{\textbf{w}}_Q\in\mathcal{M}(\overline\Omega,H\overline\Omega)$ defined as
	\begin{equation*}\label{vectormeasure}
		\int_{\overline{\Omega}}\phi(x)\cdot d\textnormal{\textbf{w}}_Q:=\int_{C([0,1],\overline{\Omega})}\left(\int_0^1\left\langle \phi(\sigma(t)),\dot{\sigma}(t)\right\rangle_H dt\right)dQ(\sigma),
	\end{equation*}
	for any continuous horizontal vector field $\phi\in C(\overline\Omega,H\overline\Omega)$. We will call this measure \textit{vector horizontal traffic intensity induced by $Q$}.
\end{deff}

The following lemma holds.
\begin{lemma}\label{7marzo24}
	Let $Q\in\mathcal{Q}_H(\mu,\nu)$, then
	\begin{enumerate}
		\item $\|\textnormal{\textbf{w}}_Q\|_{\mathcal{M}(\overline{\Omega},H\overline{\Omega})}\leq \int_{C([0,1],\overline{\Omega})}\ell(\sigma)dQ(\sigma)<+\infty;$
		\item $|\textnormal{\textbf{w}}_Q|\leq i_{Q}$ as measures;
		\item $\textnormal{div}_H\textnormal{\textbf{w}}_Q=\mu-\nu.$
	\end{enumerate}
\end{lemma}
\begin{proof}
	
(1) follows from the definition of total variation norm and of horizontal traffic plan, while (2) follows from the following inequality
\begin{equation*}
    \left|\int_{\overline{\Omega}}\phi\cdot d\textnormal{\textbf{w}}\right|\leq\int_{\overline{\Omega}}|\phi|_Hdi_Q,\quad \forall \phi\in C(\overline\Omega, H\overline\Omega).
\end{equation*}

As for (3), by definition $\textnormal{\textbf{w}}_Q$ is compactly supported, then its horizontal divergence is well defined: let $\varphi\in C^\infty(\mathbb{H}^n)$, it holds
\begin{align*}
	\int_{\overline{\Omega}}\nabla_H\varphi\cdot d\textnormal{\textbf{w}}_Q&=\int_{C([0,1],\overline{\Omega})}\left(\int_0^1\left\langle \nabla_H\varphi(\sigma(t)),\dot{\sigma}(t)\right\rangle_H dt\right)dQ(\sigma)=\\&= \int_{C([0,1],\overline{\Omega})}\left(\int_0^1\frac{d}{dt} (\varphi(\sigma(t)))dt\right)dQ(\sigma)=\\&=\int_{C([0,1],\overline{\Omega})}[\varphi(\sigma(1))-\varphi(\sigma(0))]dQ(\sigma)=\\&=\int_{\overline{\Omega}}\varphi d((e_1)_{\#}Q-(e_0)_{\#}Q)=-\int_{\overline{\Omega}}\varphi d(\mu-\nu),
\end{align*}
and the thesis follows.
\end{proof}

\begin{exa}
Following \cite[Example 4.1]{Brasco2}, we consider $x\in\mathbb{H}^n$, $y\in\mathbb{H}^n\setminus L_x,\sigma$ and $\Omega$ as in Example \ref{22marzo1}. Moreover, we consider the two discrete measures
\begin{equation*}
    \mu=\nu=\frac{\delta_x+\delta_y}{2}
\end{equation*}
and the horizontal traffic plan
\begin{equation*}
	Q:=\frac{1}{2}\delta_{\sigma}+\frac{1}{2}\delta_{\tilde\sigma}\in\mathcal{Q}_H(\mu,\nu),
\end{equation*}
where $\tilde\sigma(t):=\sigma(1-t)$. Again
\begin{equation*}
    i_Q=\mathcal{H}^1\llcorner_{\sigma([0,1])}
\end{equation*}
but
\begin{equation*}
	\textnormal{\textbf{w}}_Q\equiv0.
\end{equation*}
This shows that, in the previous Lemma, (2) is not an equality in general, because the curves of $Q$ may produce some cancellations: this is due to the fact that $\textnormal{\textbf{w}}_Q$ takes into account the orientation of the curves, while $i_Q$ does not. 
\end{exa}

\begin{lemma}\label{14marzo1}
Let $1<p<+\infty$. If 
\begin{equation*}\label{14marzo2}
	\mathcal{Q}_H^p(\mu,\nu):=\left\{Q\in\mathcal{Q}_H(\mu,\nu): i_Q\in L^p(\Omega) \right\}\not=\emptyset,
\end{equation*}
then
\begin{equation*}\label{nonemptinessbeckmann}
	\left\{\textnormal{\textbf{w}}\in L^p(\Omega,H\Omega):\textnormal{div}_H\textnormal{\textbf{w}}=\mu-\nu\right\}\not=\emptyset.
\end{equation*}
Moreover $\mu-\nu\in \left(HW^{1,q}(\Omega)\right)'$, with $q=\frac{p}{p-1}$.
\end{lemma}
\begin{proof}
Let $Q\in\mathcal{Q}_H^p(\mu,\nu)$. Lemma \ref{7marzo24} implies that $\textbf{w}_Q\in L^p(\Omega,H\Omega)$ and $\textnormal{div}_H\textbf{w}_Q=\mu-\nu$. Moreover 
\begin{multline*}\label{14marzo3}
	+\infty>\|\textnormal{\textbf{w}}_Q\|_{L^p(\Omega,H\Omega)}=\sup\left\{ \int_\Omega \left\langle\phi, \textnormal{\textbf{w}}_Q\right\rangle_H dx:\phi\in L^q(\Omega,H\Omega),\|\phi\|_{L^q(\Omega,H\Omega)}\leq1\right\}\\ \geq\sup\left\{\int_\Omega \left\langle\nabla_H\varphi, \textnormal{\textbf{w}}_Q\right\rangle_H dx:\varphi\in HW^{1,q}(\Omega),\|\varphi\|_{HW^{1,q}(\Omega)}\leq1\right\}\\=\sup\left\{\int_\Omega\varphi\ d(\mu-\nu):\varphi\in HW^{1,q}(\Omega),\|\varphi\|_{HW^{1,q}(\Omega)}\leq1\right\}=\|\mu-\nu\|_{(HW^{1,q}(\Omega))'}.
\end{multline*}
\end{proof}

As a consequence, we see that if $\mu,\nu$ are Dirac masses then $\mathcal{Q}_H^p(\mu,\nu)=\emptyset$ whenever $p\geq\frac{N}{N-1}$. Indeed, this is so because $\mu-\nu\in (HW^{1,q}(\Omega))'$ only for $q>N$, since otherwise $HW^{1,q}(\Omega)$ contains non continuous elements.

\subsection{Well-posedness of the problem}\label{characterizationofdualspacesubsection}  

In the begining of this subsection, we want to determine the image of the horizontal divergence operator $\textnormal{div}_H$ on $L^p(\Omega,H\Omega)$. To do  this, we start with the following fact.

\begin{lemma}\label{referee}
If $1<p<\infty$, and $q=\frac{p}{p-1}$, then the horizontal divergence induces a continuous linear operator $\textnormal{div}_H:L^p(\Omega,H\Omega)\to(HW^{1,q}(\Omega))'$ with norm $\leq 1$.
\end{lemma}
\begin{proof}
    Any given $\textnormal{\textbf{w}}\in L^p(\Omega,H\Omega)$  induces in a natural way a $H\mathbb{H}^n$-valued compactly supported Radon measure $\textbf{w}\in \mathcal{M}_c(\mathbb{H}^n,H\mathbb{H}^n)$. Hence the horizontal divergence $\textnormal{div}_H\textbf{w}$ is well defined, 
\begin{equation*}
	\langle \textnormal{div}_H\textnormal{\textbf{w}},\varphi\rangle
    =-\int_{\Omega}  \nabla_H\varphi \cdot d\textnormal{\textbf{w}},\quad\forall \varphi\in C^\infty(\mathbb{H}^n).
\end{equation*}
Moreover, since $\textbf{w} \in L^p(\Omega,H\Omega)$, we easily see that
\begin{equation*}
	|\langle\text{div}_H\textnormal{\textbf{w}},\varphi\rangle|\leq \|\textnormal{\textbf{w}}\|_{L^p(\Omega,H\Omega)}\|\varphi\|_{HW^{1,q}(\Omega)},\quad \forall \varphi \in HW^{1,q}(\Omega).
\end{equation*}
This allows us to extend the action of $\text{div}_H\textnormal{\textbf{w}}$ on any testing function $\varphi\in HW^{1,q}(\Omega)$, still having the above inequality. In particular, it tells us that $\text{div}_H$ has norm $\leq 1$
 as claimed. \end{proof}

The operator $\text{div}_H :L^p(\Omega,H\Omega)\to(HW^{1,q}(\Omega))'$ is not surjective. Indeed, for any $\textnormal{\textbf{w}}\in L^p(\Omega,H\Omega)$ we see that 
$$\langle \textnormal{div}_H \textnormal{\textbf{w}},1\rangle=0,$$
while there are many other elements in $(HW^{1,q}(\Omega))'$ whose action against $1$ is not $0$ (for instance, any strictly positive function in $L^p(\Omega)$). This motivates us to introduce 
 a particular closed subspace $(HW^{1,q})'_\diamond(\Omega)$ of $(HW^{1,q}(\Omega))'$,  
\begin{equation*}
	(HW^{1,q})'_\diamond(\Omega):=\left\{f\in (HW^{1,q}(\Omega))'\text{ such that }\langle f,1\rangle=0\right\}.
\end{equation*}
We endow $(HW^{1,q})'_\diamond(\Omega)$ with the environmental dual norm $\|\cdot\|_{(HW^{1,q}(\Omega))'}$. Following \cite[Section 2]{santambrogioregularity} one can characterize such space.

\begin{prop}\label{charactdualspace}
It holds that
\begin{equation*}
    (HW^{1,q})'_\diamond (\Omega)=\textnormal{div}_H\left(L^p(\Omega,H\Omega)\right).
\end{equation*}
Moreover, there is a constant $c=c(n,q,\Omega)\geq1$ such that for every $f\in (HW^{1,q})'_\diamond (\Omega)$, there is a horizontal vector field $\textnormal{\textbf{w}}\in L^p(\Omega,H\Omega)$ such that 
\begin{equation*}
	\textnormal{div}_H\textnormal{\textbf{w}}=f,
\end{equation*}
and $\|\textnormal{\textbf{w}}\|_{L^p(\Omega,H\Omega)}\leq c\|f\|_{(HW^{1,q}(\Omega))'}$.
\end{prop}

\begin{proof}
Let us consider a horizontal vector field $\textnormal{\textbf{w}}\in L^p(\Omega,H\Omega)$. We know that its weak divergence defines a distribution $\text{div}_H\textnormal{\textbf{w}}\in(HW^{1,q}(\Omega))'$ and it is trivial to prove that $\langle\text{div}_H\textnormal{\textbf{w}},1\rangle=0$. Therefore, 
$$\text{div}_H\textnormal{\textbf{w}}\in(HW^{1,q})'_\diamond(\Omega).$$
Conversely, let us consider the following minimization problem
\begin{align*}
	\min_{\varphi\in HW^{1,q}(\Omega)}\mathcal{F}(\varphi),
\end{align*}
where 
\begin{equation*}
	\varphi\mapsto \mathcal{F}(\varphi):=\frac{1}{q}\int_\Omega |\nabla_H\varphi|_H^qdx+\langle f,\varphi\rangle.
\end{equation*}
This problem admits at least one solution. Indeed, we can restrict the minimization to the subspace
\begin{equation}\label{7dicembre2}
	HW_\diamond^{1,q}(\Omega):=\left\{\varphi\in HW^{1,q}(\Omega):\int_\Omega\varphi(x)dx=0\right\},
\end{equation}
which is a convex subset and closed w.r.t. weak topology of $HW^{1,q}(\Omega)$.	The functional $\mathcal{F}$ is lower semicontinuous w.r.t. the weak topology of $HW^{1,q}(\Omega)$ and it is coercive on $HW_\diamond^{1,q}(\Omega)$: indeed, if $\varphi\in HW_\diamond^{1,q}(\Omega)$, then from \eqref{poincaréineq} it follows the existence of a constant $c=c(n,q,\Omega)>0$ such that
\begin{equation*}
\aligned
    |\mathcal{F}(\varphi)|&\geq\frac{1}{q}\|\nabla_H\varphi\|^q_{L^q(\Omega,H\Omega)}-(c^{1/q}+1)\|f\|_{(HW^{1,q}(\Omega))'}\|\nabla_H\varphi\|_{L^q(\Omega,H\Omega)}\\&=\|\nabla_H\varphi\|_{L^q(\Omega,H\Omega)}\left(\frac{1}{q}\|\nabla_H\varphi\|_{L^q(\Omega,H\Omega)}^{q-1}-(c^{1/q}+1)\|f\|_{(HW^{1,q}(\Omega))'}\right).
\endaligned
\end{equation*}
Hence, the existence of minimizers follows again from \eqref{poincaréineq}.
	
Computing Euler Lagrange equations, a solution $\varphi$ of the problem above turns out to satisfy 
\begin{equation*}
	-\int_\Omega \langle|\nabla_H\varphi|^{q-2}\nabla_H\varphi,\nabla_H\psi \rangle_Hdx=\langle f,\psi\rangle, \quad \forall \psi\in HW^{1,q}.
\end{equation*}
This means that there exists $\textnormal{\textbf{w}}=|\nabla_H\varphi|^{q-2}\nabla_H\varphi\in L^p(\Omega,H\Omega)$ such that $\text{div}_H\textnormal{\textbf{w}}=f$. Moreover, testing $f$ against $\varphi$ we get
\begin{equation*}
\aligned
    \|\textnormal{\textbf{w}}\|^p_{L^p(\Omega,H\Omega)}&=\int_\Omega |\textnormal{\textbf{w}}|_H^pdx=\int_\Omega|\nabla_H\varphi|_H^qdx =-\langle f,\varphi\rangle\\&\leq \|f\|_{(HW^{1,q}(\Omega))'}\|\varphi\|_{HW^{1,q}(\Omega)} \leq (c^{1/q}+1)\|f\|_{(HW^{1,q}(\Omega))'}\|\nabla_H\varphi\|_{L^q(\Omega,H\Omega)}\\&=(c^{1/q}+1)\|f\|_{(HW^{1,q}(\Omega))'}\|\textnormal{\textbf{w}}\|_{L^p(\Omega,H\Omega)}^{p-1},
\endaligned
\end{equation*}
as desired.
\end{proof}

Let now consider a function $\mathcal{G}:\overline{\Omega}\times\mathbb{R}^{2n}\to\mathbb{R}$ such that
\begin{equation*}
	x\mapsto\mathcal{G}(x,w)\textnormal{ is Lebesgue-measurable for any }w\in\mathbb{R}^{2n}.
\end{equation*}
Let us suppose in addition that $\mathcal{G}(x,0)=0$ for any $x\in\overline{\Omega}$,
\begin{equation}\label{hypconvexity}
	w\mapsto\mathcal{G}(x,w) \text{ is convex for any }x\in\overline\Omega
\end{equation}
and it satisfies
\begin{equation}\label{pgrowthmathcalG}
	\frac{a}{p}|w|^p-h_0(x)\leq\mathcal{G}(x,w)\leq \frac{b}{p}|w|^p+h_1(x), \textnormal{ for any }(x,w)\in\overline{\Omega}\times\mathbb{R}^{2n}
\end{equation}
for some $h_0,h_1\in L^1(\Omega)$, $h_0,h_1\geq0$, some constants $a,b>0$ and $p\in(1,\infty)$.

Following \cite{santambrogioregularity}, we are interested in the Beckmann-type problem
\begin{equation}\label{generalbeckmann}\tag{$\mathcal{B}$}
	\inf_{\textnormal{\textbf{w}}\in L^p(\Omega,H\Omega) }\left\{\int_\Omega\mathcal{G}(x,\textnormal{\textbf{w}}(x))dx:\textnormal{div}_H\textnormal{\textbf{w}}=f \right\}.
\end{equation}
Since $\textnormal{div}_H\textnormal{\textbf{w}}\in (HW^{1,q})'_\diamond(\Omega)$ for any $\textnormal{\textbf{w}}\in L^p(\Omega,H\Omega)$, the natural functional setting for the above problem is $f \in(HW^{1,q})'_\diamond(\Omega)$.

\begin{teo}\label{existenceofbeckmin}
If $f\in (HW^{1,q})'_\diamond (\Omega)$, then   \eqref{generalbeckmann} admits a minimizer. If $\mathcal{G}$ is strictly convex, then the minimizer is unique.
\end{teo}
\begin{proof}
Let $f\in (HW^{1,q})'_\diamond (\Omega)$, then Proposition \ref{charactdualspace} implies that  there exists at least one vector field $\textnormal{\textbf{w}}\in L^p(\Omega,H\Omega)$ such that $\textnormal{div}_H\textnormal{\textbf{w}}=f$: hence \eqref{pgrowthmathcalG} implies that $\eqref{generalbeckmann}<+\infty$. Let $(\textnormal{\textbf{w}}_n)_{n\in\mathbb{N}}\subseteq L^p(\Omega,H\Omega)$ be a minimizing sequence. Since we are dealing with a minimizing sequence, from \eqref{pgrowthmathcalG} it follows that this sequence is bounded in $L^p(\Omega,H\Omega)$: hence, up to a subsequence, it is weakly convergent to some $\tilde{\textnormal{\textbf{w}}}\in L^p(\Omega,H\Omega)$. This vector field is still admissible: indeed, by the weak convergence it follows that
\begin{equation*}
    -\int_\Omega \left\langle \tilde{\textnormal{\textbf{w}}}, \nabla_H\varphi\right\rangle_Hdx=-\lim_{n\to\infty}\int_\Omega \left\langle\textnormal{\textbf{w}}_n, \nabla_H\varphi\right\rangle_Hdx=\langle f,\varphi\rangle,\quad \forall\varphi\in HW^{1,q}(\Omega).
\end{equation*}
Since $\mathcal{G}$ is convex in the second variable, the integral functional is convex, as well, and lower semicontinuous; this implies the lower semicontinuity of the integral functional w.r.t. the weak convergence in $L^p(\Omega,H\Omega)$, hence
\begin{equation*}
    \int_\Omega\mathcal{G}(x,\tilde{\textnormal{\textbf{w}}}(x))dx\leq\liminf_{n\to+\infty}\int_\Omega\mathcal{G}(x,\textnormal{\textbf{w}}_n(x))dx=\eqref{generalbeckmann}<+\infty. 
\end{equation*}
This proves that $\tilde{\textnormal{\textbf{w}}}$ is a minimum.
\end{proof}

\subsection{Dual formulation}
The aim of this subsection is to show that the Beckmann-type problem \eqref{generalbeckmann} admits the dual formulation
\begin{equation}\label{generaldualbeckmann}\tag{$\mathcal{D}$}
	\sup _{\varphi\in HW^{1,q}(\Omega)}\left\{-\left\langle f,\varphi\right\rangle-\int_\Omega \mathcal{G}^*(x,\nabla_H\varphi(x))dx\right\},
\end{equation}
where $f\in (HW^{1,q})'_\diamond(\Omega)$, $\mathcal{G}^*$ is the Legendre transform of the map $w\mapsto\mathcal{G}(x,w)$ (see Theorem \ref{convexdualitythm}), and again we are identifying $\nabla_H\varphi=\sum_{i=1}^{2n}X_i\varphi X_i$ with $(X_1\varphi,\ldots,X_{2n}\varphi):\Omega\to\mathbb{R}^{2n}$.

By definition the function
\begin{equation*}
	x\mapsto\mathcal{G}^*(x,z)\textnormal{ is Lebesgue-measurable for any }z\in\mathbb{R}^{2n}
\end{equation*}
and
\begin{equation}\label{hypconvexitytransform}
	\mathbb{R}^{2n}\ni z\mapsto \mathcal{G}^*(x,z)\text{ is convex, for any }x\in\overline{\Omega}.
\end{equation}
The condition \eqref{pgrowthmathcalG} implies the following $q$-growth condition on $\mathcal{G}^*$
\begin{equation}\label{qgrowthmathcalG*}
	\frac{1}{qb^{q-1}}|z|^q-h_1(x)\leq\mathcal{G}^*(x,z)\leq\frac{1}{qa^{q-1}}|z|^q+h_0(x),\textnormal{ for any }(x,z)\in \overline{\Omega}\times\mathbb{R}^{2n}
\end{equation}
where $q=\frac{p}{p-1}$.

We follow the proof developed in \cite{Brasco}, which is based on the following classical theorem of convex analysis.

\begin{teo}[Convex duality]\cite[Proposition 5]{ekeland2012convexity}\label{convexdualitythm}
Let $\mathcal{F}:Y\to\mathbb{R}$ a convex and lower semicontinuous functional on a reflexive Banach space $Y$. Let $X$ be another reflexive Banach space and $A:X\to Y$ a bounded linear operator, with adjoint operator $A':Y'\to X'$. Then we have
\begin{equation}\label{convexduality}
	\sup_{x\in X}\langle x',x\rangle-\mathcal{F}(Ax)=\inf_{y'\in Y'}\left\{\mathcal{F}^*(y'): A'y'=x'\right\},\quad x'\in X',
\end{equation}
where $\mathcal{F}^*:Y'\to\mathbb{R}$ denotes the Legendre-Fenchel transform of $\mathcal{F}$.
		
Moreover, if the supremum in \eqref{convexduality} is attained at some $x_0\in X$, then infimum in \eqref{convexduality} is attained at some $y_0'\in Y'$ such that 
\begin{equation}\label{primaldual}
	y_0'\in \partial \mathcal{F}(Ax_0),
\end{equation}
where $\partial\mathcal{F}$ denotes the subdifferential
\begin{equation*}
	\partial\mathcal{F}(y):=\left\{y'\in Y':\mathcal{F}(y)+\mathcal{F}^*(y')=\langle y',y\rangle\right\},\quad y\in Y.
\end{equation*}
\end{teo}

We are now ready to show that $\eqref{generalbeckmann}=\eqref{generaldualbeckmann}$.

\begin{teo}\label{dualteo}
If $f\in (HW^{1,q})'_\diamond (\Omega)$, then the problem \eqref{generaldualbeckmann} admits a solution and
\begin{equation*}
	\min_{\textnormal{\textbf{w}}\in L^p(\Omega,H\Omega) }\left\{\int_\Omega\mathcal{G}(x,\textnormal{\textbf{w}}(x))dx:\textnormal{div}_H\textnormal{\textbf{w}}=f \right\}=\max_{\varphi\in HW^{1,q}(\Omega)}\left\{-\int_\Omega \mathcal{G}^*(x,\nabla_H \varphi(x))dx-\langle f,\varphi\rangle\right\}.
\end{equation*}
Moreover, if $\textnormal{\textbf{w}}_0\in L^p(\Omega,H\Omega)$ in a minimizer for \eqref{generalbeckmann} and $\varphi_0\in HW^{1,q}(\Omega)$ is a maximizer for \eqref{generaldualbeckmann}, then we have the following primal-dual optimality condition
\begin{equation}\label{reciprocityformulabeckmannanddual}
	\textnormal{\textbf{w}}_0\in\partial\mathcal{G}^*(x,\nabla_H\varphi_0(x)) \text{ a.e. in }\Omega.
\end{equation}
\end{teo}

\begin{proof}
First we observe that the existence of a minimizer for the problem \eqref{generalbeckmann} follows from Theorem \ref{existenceofbeckmin}. Thanks to the Direct method in calculus of variations also \eqref{generaldualbeckmann} admits at least a solution belonging the space $HW^{1,q}_\diamond(\Omega)$ of Sobolev functions with zero mean. Indeed, \eqref{hypconvexitytransform} implies that the functional $\mathcal{\tilde{F}}:HW^{1,q}(\Omega)\to\mathbb{R}$
\begin{equation*}
	\mathcal{\tilde{F}}(\varphi)=\int_\Omega\mathcal{G}^*(x,\nabla_H\varphi(x))dx+\langle f,\varphi\rangle.
\end{equation*}
is lower semicontinuous with respect to the weak topology of $HW^{1,q}(\Omega)$; moreover it is coercive on $HW^{1,q}_{\diamond}(\Omega)$: let $\varphi\in HW^{1,q}_{\diamond}(\Omega)$ from \eqref{qgrowthmathcalG*} and \eqref{poincaréineq}  
\begin{multline*}
	\mathcal{\tilde{F}}(\varphi)\geq \frac{1}{qb^{q-1}}\|\nabla_H\varphi\|_{L^q(\Omega,H\Omega)}^{q}-\|f\|_{(HW^{1,q}(\Omega))'}\|\varphi\|_{HW^{1,q}(\Omega)}-d\\ \geq \frac{1}{qb^{q-1}}\|\nabla_H\varphi\|_{L^q(\Omega,H\Omega)}^{q}-(c^{1/q}+1)\|f\|_{(HW^{1,q}(\Omega))'}\|\nabla_H\varphi\|_{L^q(\Omega,H\Omega)}-d\\
	\geq\|\nabla_H\varphi\|_{L^q(\Omega,H\Omega)}\left(\frac{1}{qb^{q-1}}\|\nabla_H\varphi\|_{L^q(\Omega,H\Omega)}^{q-1}-(c^{1/q}+1)\|f\|_{(HW^{1,q}(\Omega))'}\right)-d.
\end{multline*}
where $d>0$ and $c=c(n,q,\Omega)$ is the constant appearing in \eqref{poincaréineq}.
The existence of solutions follows from again from \eqref{poincaréineq}.
		
Let now denote by $X=HW^{1,q}(\Omega)$, $Y=L^q(\Omega,H\Omega)$ and let us consider the operator
$$A:X\to Y,\quad A(\varphi)=\nabla_H\varphi,\quad \forall \varphi\in X.$$
This operator is linear; moreover it is bounded since 
$$\|A(\varphi)\|_Y=\|\nabla_H\varphi\|_{L^q(\Omega,H\Omega)}\leq \|\varphi\|_{HW^{1,q}(\Omega)}.$$
We denote by $\mathcal{F}:Y\to \mathbb{R}$ the functional
$$\mathcal{F}(\phi):=\int_\Omega \mathcal{G}^*(x,\phi(x))dx.$$
By the convexity and the lower semicontinuity of $\mathcal{G}$ it follows that $\mathcal{F}^*:L^p(\Omega,H\Omega)\to \mathbb{R}$
$$\mathcal{F}^*(\textnormal{\textbf{w}})=\int_\Omega \mathcal{G}^{**}(x,\textnormal{\textbf{w}}(x))dx=\int_\Omega \mathcal{G}(x,\textnormal{\textbf{w}}(x))dx.$$
Moreover, it is a straightforward consequence of the definition of the continuous linear operator
\begin{equation*}
    \textnormal{div}_H:L^p\left(\Omega,H\Omega\right)\longrightarrow (HW^{1,q}(\Omega))' 
\end{equation*}
that $A'=-\textnormal{div}_H$.
		
The rest of the thesis follows from the primal-dual optimality condition \eqref{primaldual}.
		
\end{proof}

\begin{Remark}
In \cite{santambrogioregularity} the author provides a different proof of the equality 
\begin{equation*}
	\eqref{generalbeckmann}=\eqref{generaldualbeckmann},
\end{equation*}
which works also in the Heisenberg group. For the sake of completeness we just highlight the main steps.  
Let $\mathcal{F}:(HW^{1,q}(\Omega))'\to\mathbb{R}$ be the functional
\begin{equation*}
	\mathcal{F}(k):=\min_{\textnormal{\textbf{w}}\in L^p(\Omega,H\Omega)}\left\{ \int_\Omega \mathcal{G}(x,\textnormal{\textbf{w}}(x))dx: \textnormal{div}_H\textnormal{\textbf{w}}=f+k\right\}. 
\end{equation*}
Lemma \ref{charactdualspace} and \eqref{pgrowthmathcalG} imply that the functional $\mathcal{F}(k)$ is well-defined and real-valued if and only if $k\in (HW^{1,q})'_\diamond(\Omega)$.
A computation implies that its Legendre transform $\mathcal{F}^*:HW^{1,q}(\Omega)\to\mathbb{R}$ is precisely 
\begin{equation*}
	\mathcal{F}^*(\varphi)=-\langle f,\varphi\rangle+\int_\Omega \mathcal{G}^*(x,-\nabla_H\varphi(x))dx.
\end{equation*}
By definition $\mathcal{F}^{**}(0)=\sup_{\varphi\in HW^{1,q}(\Omega)}\left\{-\mathcal{F}^*(\varphi)\right\}$; again, since one can restrict the minimization to $HW^{1,q}_\diamond(\Omega)$ hence $\sup-\mathcal{F}^*<+\infty$. By taking the sup on $-\varphi$ instead of $\varphi$ we also have
\begin{align*}
	\mathcal{F}^{**}(0)=\sup_{\varphi\in HW^{1,q}(\Omega)}\left\{-\langle f,\varphi\rangle-\int_\Omega \mathcal{G}^*(x,\nabla_H\varphi(x))dx\right\}.
\end{align*}
One can prove that $\mathcal{F}$ is convex and l.s.c., then $\mathcal{F}^{**}(0)=\mathcal{F}(0)$ and the thesis follows.
\end{Remark}

Let us remark that if the function $w \mapsto \mathcal{G}(x, w)$ is strictly convex, then the map $z \mapsto \mathcal{G}^*(x, z)$ is $C^1$ (see \cite[Proposition 4.26]{santambrogiolibro2}). Therefore, its subdifferential reduces to the singleton $D\mathcal{G}^*$ and, according to \eqref{reciprocityformulabeckmannanddual} in Theorem \ref{dualteo}, the unique solution to \eqref{generalbeckmann} can be written as
\begin{equation*}
\textnormal{\textbf{w}}_0 = D\mathcal{G}^*(x, \nabla_H \varphi_0(x)) \text{ a.e. in } \Omega,
\end{equation*}
where $\varphi_0$ is a weak solution to the Euler-Lagrange equation
\begin{equation}\label{beckeul}
\textnormal{div}_H \left( D\mathcal{G}^*(x, \nabla_H \varphi(x)) \right) = f.
\end{equation} 

\begin{exa}
If $\mathcal{G}(w)=\frac{1}{p}|w|^p$, then $\mathcal{G}^*(z)=\frac{1}{q}|z|^q$ and \eqref{beckeul} becomes the degenerate $q$-Laplace equation 
\begin{equation}\label{12agosto}
	\textnormal{div}_H\left(|\nabla_H\varphi|_H^{q-2}\nabla_H\varphi\right)=f,\quad \text{in }\Omega.
\end{equation}
Based on the results in \cite{Capogna} and \cite[Theorem 3.1]{zhong2017regularity}, it can be shown that if $f, X_{2n+1}f \in L_{loc}^2(\Omega)$, and $\varphi_0 \in HW^{1,q}(\Omega)$ is a weak solution to \eqref{12agosto}, then
\begin{equation*}
    F := |\nabla_H \varphi_0|^{\frac{q-2}{2}} \nabla_H \varphi_0 \in HW_{loc}^{1,2}(\Omega, H\Omega).
\end{equation*}
Consequently, following \cite{Brasco}, the solution $\textnormal{\textbf{w}}_0$ to \eqref{generalbeckmann} satisfies
\begin{equation*}
    \textnormal{\textbf{w}}_0=D\mathcal{G}^*(\nabla_H\varphi_0)=|F|^{\frac{q-2}{q}}F\in HW_{loc}^{1,r}(\Omega,H\Omega),
\end{equation*}
where 
\begin{equation*}
    r:=
    \begin{cases}
    2&\textnormal{if }q=2,\\
    \frac{Nq}{(N-1)q+2-N}&\textnormal{if }2<q<4.
    \end{cases}
\end{equation*}
On the other hand, the optimal regularity for solutions is the H\"older continuity of the horizontal derivatives: for $2\leq q<\infty$, it has been established in \cite{mukherjee2019regularity}, with $f\in L^p_{loc}(\Omega)$ and $p>N$.

See also \cite{circelliregularity} for Lipschitz regularity for solutions to \eqref{beckeul}, with a different choice of $\mathcal{G}$ and $f=0$, and \cite{circelli2025gradientestimatesorthotropicnonlinear} for its parabolic counterpart.
\end{exa}

\subsection{Lagrangian formulation}

The \textit{congested optimal transport problem}, adapted to the Heisenberg group in \cite{circelli2023transport}, consists of the following minimization problem:
\begin{equation}\label{lepbme1chap4}\tag{$\mathcal{W}$}
	\inf_{Q\in\mathcal{Q}^p_H(\mu,\nu)}\int_\Omega G(x,i_Q(x)) \, dx,
\end{equation}
where the function $G:\overline{\Omega}\times\mathbb{R}\to\mathbb{R}$ satisfies:
\begin{equation}\label{hyp-1}
	x \mapsto G(x,i) \text{ is Lebesgue-measurable for any } i \in \mathbb{R},
\end{equation}
and
\begin{equation}\label{hyp0}
	G(x,0) = 0, \text{ for any } x \in \overline{\Omega}.
\end{equation}
Additionally, we assume:
\begin{equation}\label{convexityG}
	i \mapsto G(x,i) \text{ is convex and non-decreasing for any } x \in \overline{\Omega},
\end{equation}
and it exhibits $p$-growth in the second variable:
\begin{equation}\label{Gbehaviour}
	\frac{a}{p}i^p - h_0(x) \leq G(x,i) \leq \frac{b}{p} i^p + h_1(x), \text{ for any } (x,i) \in \overline{\Omega} \times \mathbb{R},
\end{equation}
for some $p \in (1,\infty)$, where $a, b > 0$ and $h_0, h_1 \in L^p(\Omega, \mathbb{R}_+)$. It turns out that solutions to \eqref{lepbme1chap4} can, in some cases, exhibit the nature of Wardrop equilibria. For further details, we refer the interested reader to \cite{circelli2023transport}.

In this subsection, we investigate the equivalence between \eqref{lepbme1chap4} and \eqref{generalbeckmann}, where $f = \mu - \nu$, and the function $\mathcal{G}(x,\cdot) = G(x,|\cdot|_H)$ for any $x \in \overline{\Omega}$:
\begin{equation*}
    \inf_{\textnormal{\textbf{w}}\in L^p(\Omega,H\Omega)} \left\{ \int_\Omega G(x,|\textnormal{\textbf{w}}(x)|_H) \, dx : \textnormal{div}_H \textnormal{\textbf{w}} = \mu - \nu \right\}.
\end{equation*}
Moreover, we will demonstrate how to pass from a solution of one problem to a solution of the other, and vice versa.

First, we state the following result.

\begin{prop}\label{21marzo251}
If $\mathcal{Q}_H^p(\mu,\nu) \neq \emptyset$, then $\mu - \nu \in (HW^{1,q})'_\diamond(\Omega)$. More precisely, for every $Q \in \mathcal{Q}_H^p(\mu,\nu)$, $\textnormal{\textbf{w}}_Q$ is admissible for \eqref{generalbeckmann} and $|\textnormal{\textbf{w}}_Q|_H \leq i_Q$. Therefore,
\begin{equation}\label{4marzo242}
    \eqref{generalbeckmann} \leq \eqref{lepbme1chap4}.
\end{equation}
\end{prop}

\begin{proof}
The result follows from Lemma \ref{14marzo1}, Proposition \ref{charactdualspace} and the monotonicity of $G$.
\end{proof}

Hence, the equivalence between the two problems is reduced to the converse inequality.

To prove this, we follow \cite[Theorem 2.1]{brasco2013congested} and \cite{santambrogio2014dacorogna}, which are based on the Dacorogna-Moser scheme combined with a regularization procedure. In the Euclidean setting, this regularization procedure is unnecessary if $\mu$ and $\nu$ are sufficiently smooth, because every solution $\textnormal{\textbf{w}}$ to \eqref{generalbeckmann} belongs to some Sobolev space. Hence, one can define the flow in the weaker sense of DiPerna-Lions, as seen in \cite{Brasco}. The absence of a sub-Riemannian DiPerna-Lions theory, as discussed in \cite{Ambrosio}, is one of the reasons to force a regularization in our setting. 

\begin{teo}\label{9dicembre2}
If $\mu - \nu \in (HW^{1,q})'_\diamond(\Omega)$, then $\mathcal{Q}_H^p(\mu,\nu) \neq \emptyset$. More precisely, for every $\textnormal{\textbf{w}}$ admissible for \eqref{generalbeckmann}, there is $Q_{\textnormal{\textbf{w}}}\in\mathcal{Q}_H^p(\mu,\nu)$ such that $i_{Q_{\textnormal{\textbf{w}}}}\leq|\textnormal{\textbf{w}}|_H$. 
Therefore,
\begin{equation*}
    \eqref{generalbeckmann} \geq \eqref{lepbme1chap4}.
\end{equation*}
\end{teo}

\begin{proof}
Let $\textnormal{\textbf{w}}=\sum_{j=1}^{2n}\textnormal{\textbf{w}}_jX_j\in L^p(\Omega,H\Omega)$ be a solution to \eqref{generalbeckmann}. We extend it by $0$ outside $\Omega$ and we consider the horizontal vector field $\textnormal{\textbf{w}}^\varepsilon=\sum_{j=1}^{2n}\textnormal{\textbf{w}}^\varepsilon_jX_j$, where
\begin{equation*}
	\textnormal{\textbf{w}}^\varepsilon_j(x):=\rho_\varepsilon\ast\textnormal{\textbf{w}}_j(x),\quad \forall j=1,\ldots,2n,
\end{equation*}
$\rho_\varepsilon$ is a mollifier for the group structure of $\mathbb{H}^n$ and $\Omega_\varepsilon:=B(0,\varepsilon)\cdot\Omega$. We notice that, for any $j=1,\ldots,2n$, $\textnormal{\textbf{w}}^\varepsilon_j\in C^\infty_c(\mathbb{H}^n)$ with support in $\Omega_\varepsilon$. A simple computation implies that 
$$\textnormal{div}_H\textnormal{\textbf{w}}^\varepsilon=\mu^\varepsilon-\nu^\varepsilon,$$
where  
$$\mu^\varepsilon=\rho_\varepsilon\ast\mu+\varepsilon,\quad \nu^\varepsilon=\rho_\varepsilon\ast\nu+\varepsilon\in C^\infty,$$
and we are supposing that both $\mu$ and $\nu$ are extended by $0$ outside $\Omega$. Let us introduce the non-autonomous horizontal vector field 
$$\hat{\textnormal{\textbf{w}}}^\varepsilon(t,x):=\frac{\textnormal{\textbf{w}}^\varepsilon(x)}{\overline{\mu}^\varepsilon(t,x)},\quad \forall (t,x)\in[0,1]\times\overline{\Omega}_\varepsilon,$$
where  
\begin{equation}
    \overline{\mu}^\varepsilon(t,x):=(1-t)\mu^\varepsilon(x)+t\nu^\varepsilon(x)>\varepsilon>0,\quad (t,x)\in[0,1],\times\overline{\Omega}_\varepsilon.
\end{equation}
Let us notice that $\hat{\textnormal{\textbf{w}}}^\varepsilon(t,\cdot)\in C^\infty_c(\mathbb{H}^n,H\mathbb{H}^n)$ with support in $\Omega_\varepsilon$, for any $t\in[0,1]$. Moreover, for any $\epsilon>0$ it holds that
\begin{equation*}
	\mu^\varepsilon-\nu^\varepsilon=\textnormal{div}_H\textnormal{\textbf{w}}^\varepsilon=\textnormal{div}_\epsilon\textnormal{\textbf{w}}^\varepsilon,
\end{equation*}
where $\textnormal{div}_\epsilon$ denotes the Riemannian $\epsilon$-divergence in $\mathbb{H}^n\equiv\mathbb{R}^{2n+1}$, namely
$$\int\varphi\textnormal{div}_\epsilon\textnormal{\textbf{w}}^\varepsilon dx=-\int g_\epsilon\left(\textnormal{\textbf{w}}^\varepsilon,\nabla_\epsilon\varphi\right)dx$$
for any $\varphi\in C^\infty(\mathbb{H}^n)$, see also Section \ref{riemannianapproxsubsection}. In particular, by noticing that $\hat{\textnormal{\textbf{w}}}^\varepsilon(t,\cdot)\in C^\infty_c(\mathbb{R}^{2n+1},\mathbb{R}^{2n+1})$, for any $t\in[0,1]$, the curve $\overline{\mu}^\varepsilon(t,\cdot)$ satisfies the following initial value problem for the Riemannian continuity equation
\begin{equation}\label{11marzo}
	\begin{cases}\partial_t\lambda+\text{div}_\epsilon(\hat{\textnormal{\textbf{w}}}^\varepsilon\lambda)=0,\\
    \lambda(0,\cdot)=\mu^\varepsilon.
    \end{cases}
\end{equation}
Since $\hat{\textnormal{\textbf{w}}}^\varepsilon(t,\cdot)$ is smooth for any $t\in[0,1]$, the unique solution to \eqref{11marzo} is given by $\left(\Psi^\varepsilon(t,\cdot)\right)_{\#}\mu^\varepsilon$, where  $$\Psi^\varepsilon:[0,1]\times\mathbb{R}^{2n+1}\to\mathbb{R}^{2n+1}$$ is the flow of the vector field $\hat{\textbf{w}}^\varepsilon$ in the Riemannian sense,
\begin{equation*}
\begin{cases}
    \frac{d}{dt}\Psi^\varepsilon(t,x)=\hat{\textnormal{\textbf{w}}}^\varepsilon(t,\Psi_\varepsilon(t,x)),\\
    \Psi^\varepsilon(0,x)=x.
\end{cases}
\end{equation*}
Let us remark that, since $\hat{\textnormal{\textbf{w}}}^\varepsilon(t,\cdot)$ has compact support in $\Omega_\varepsilon$, for any $t\in[0,1]$, then  $\Psi^\varepsilon(t,x)\in\Omega_\varepsilon$ if $x\in\Omega_\varepsilon$, and $\Psi^\varepsilon(t,x)=x$ otherwise.

Hence
\begin{equation*}
	\left(\Psi^\varepsilon(t,\cdot)\right)_{\#}\mu^\varepsilon=\overline{\mu}^\varepsilon(t,\cdot),\quad \forall t\in[0,1].
\end{equation*}
If we denote by 
\begin{equation*}
	\Phi^\varepsilon: \overline{\Omega}_\varepsilon\to C([0,1],\overline{\Omega}_\varepsilon),
\end{equation*}
the trajectory map, i.e. $\Phi^\varepsilon(x):=\Psi^\varepsilon(\cdot,x)$, then the measure
\begin{equation*}
	Q_\varepsilon:=\left(\Phi^\varepsilon \right)_{\#}\mu^\varepsilon\in \mathcal{M}_+(C([0,1],\overline{\Omega}_\varepsilon)),
\end{equation*}
satisfies
\begin{equation*}
	\left(e_t\right)_{\#}Q_\varepsilon=\left(e_t\circ\Phi^\varepsilon \right)_\#\mu^\varepsilon=\left(\Psi^\varepsilon(t,\cdot)\right)_{\#}\mu^\varepsilon=\overline{\mu}^\varepsilon(t,\cdot),\quad \forall t\in[0,1];
\end{equation*}
in particular,
\begin{equation}\label{approx1}
	\left(e_0\right)_{\#}Q_\varepsilon=\mu^\varepsilon\quad\left(e_1\right)_{\#}Q_\varepsilon=\nu^\varepsilon.
\end{equation}
By construction $Q_\varepsilon$ is supported on the integral curves of the horizontal vector field $\hat{\textbf{w}}^\varepsilon$ and the definition of $Q_\varepsilon=\Phi^\varepsilon_\#\mu^\varepsilon$ and a change of variables imply
\begin{multline*}
	\int_{C([0,1],\overline{\Omega}_\varepsilon)}\ell_{SR}(\sigma)dQ_\varepsilon(\sigma)=\int_{\overline{\Omega}_\varepsilon}\left(\int_0^1|\hat{\textbf{w}}^\varepsilon(t,\Psi^\varepsilon(t,x))|_Hdt\right)\mu^\varepsilon(x)dx\\=\int_0^1\left(\int_{\overline{\Omega}_\varepsilon}|\hat{\textnormal{\textbf{w}}}^\varepsilon(t,y)|_Hd\Psi^\varepsilon(t,\cdot)_\#\mu^\varepsilon(y)\right)dt=\int_0^1\left(\int_{\overline{\Omega}_\varepsilon}|\hat{\textbf{w}}^\varepsilon(t,y)|_H\overline{\mu}^\varepsilon(t,y)dy\right)dt\\=\int_{\overline{\Omega}_\varepsilon}|\textbf{w}^\varepsilon(y)|_Hdy\leq\|\textnormal{\textbf{w}}\|_{L^1(\Omega,H\Omega)}\leq \mathcal{L}^{2n+1}(\Omega)^{1/q}\|\textnormal{\textbf{w}}\|_{L^p(\Omega,H\Omega)}<+\infty,
\end{multline*}
because $\textbf{w}\in L^p(\Omega,H\Omega)$.
	
Hence, one can define the horizontal traffic intensity $i_{Q_\varepsilon}$ associated with $Q_\varepsilon$ as
\begin{equation*}
	\int_{\overline{\Omega}_\varepsilon}\varphi(x)di_{Q_\varepsilon}(x):=\int_{C([0,1],\overline{\Omega}_\varepsilon)}\left(\int_0^1\varphi(\sigma(t))|\dot\sigma(t)|_Hdt\right)dQ_\varepsilon(\sigma),
\end{equation*}
for any $\varphi\in C(\overline{\Omega}_\varepsilon)$. Moreover, the definition of $Q_\varepsilon$ and a change of variables imply 
\begin{multline}\label{16marzo}
	\int_{\overline{\Omega}_\varepsilon}\varphi(x)di_{Q_\varepsilon}(x)=\int_{\overline{\Omega}_\varepsilon}\left(\int_0^1 \varphi(\Psi^\varepsilon(t,x))|\hat{\textnormal{\textbf{w}}}^\varepsilon(t,\Psi^\varepsilon(t,x))|_Hdt\right)\mu^\varepsilon(x)dx\\=\int_0^1\left(\int_{\overline{\Omega}_\varepsilon}\varphi(y)|\hat{\textnormal{\textbf{w}}}^\varepsilon(t,y)|_H\overline{\mu}^\varepsilon(t,y)dy\right)dt
	=\int_{\overline{\Omega}_\varepsilon}\varphi(y)|\textnormal{\textbf{w}}^\varepsilon(y)|_Hdy,
\end{multline}
for any $\varphi\in C(\overline{\Omega}_\varepsilon)$; hence $i_{Q_\varepsilon}=|\textnormal{\textbf{w}}^\varepsilon|_H\in C_c^\infty(\mathbb{H}^n)$, with support in $\Omega_\varepsilon$, and again
\begin{equation}\label{5dicembre}
	\|i_{Q_\varepsilon}\|_{L^1(\Omega_\varepsilon)}= \|\textnormal{\textbf{w}}^\varepsilon\|_{L^1(\Omega_\varepsilon,H\Omega_\varepsilon)}\leq \mathcal{L}^{2n+1}(\Omega)^{1/q}\|\textnormal{\textbf{w}}\|_{L^p(\Omega,H\Omega)}<+\infty. 
\end{equation}
Analogously we may define and compute the vector traffic intensity $\textbf{w}_{Q_\varepsilon}$ associated with $Q_\varepsilon$: given $\phi\in C(\overline{\Omega}_\varepsilon,H\overline{\Omega}_\varepsilon)$
\begin{equation*}
	\int_{\overline{\Omega}_\varepsilon}\phi(x)\cdot d\textbf{w}_{Q_\varepsilon}(x):=\int_{\overline{\Omega}_\varepsilon}\left(\int_0^1 \langle\phi(\Psi^\varepsilon(t,x)),\hat{\textnormal{\textbf{w}}}^\varepsilon(t,\Psi^\varepsilon(t,x))\rangle_Hdt\right)\mu^\varepsilon(x)dx=\int_{\overline{\Omega}_\varepsilon}\langle\phi(y),\textnormal{\textbf{w}}^\varepsilon(y)\rangle_Hdy,
\end{equation*}
hence $\textnormal{\textbf{w}}_{Q_\varepsilon}=\textnormal{\textbf{w}}^\varepsilon\in C_c(\mathbb{H}^n,H\mathbb{H}^n)$, with support in $\Omega_\varepsilon$.

For any $0<\varepsilon<1$, we denote by $\tilde{Q}_\varepsilon:=Q_\varepsilon/(1+\varepsilon)$ and we fix a compact set $\Omega'\subset\mathbb{H}^n$, such that $\Omega\subset\Omega_\varepsilon\subset\Omega'$. Moreover, it holds that
\begin{equation*}
	\tilde{Q}_\varepsilon(C([0,1],\Omega'))=\frac{Q_\varepsilon(C([0,1],\Omega'))}{1+\varepsilon}=\frac{1}{1+\varepsilon}\int_{\Omega'}\mu^\varepsilon(x)dx=\frac{1}{1+\varepsilon}\int_{\Omega_\varepsilon}\mu^\varepsilon(x)dx=1,
\end{equation*}
therefore $\tilde{Q}_\varepsilon\subseteq\mathcal{P}\left(C([0,1],\Omega')\right)$. Since both $i_{\tilde{Q}_\varepsilon}$ and $\textnormal{\textbf{w}}_{\tilde{Q}_\varepsilon}$ are invariant by reparametrization, we may suppose that $\tilde{Q}_\varepsilon$ is supported on curves parametrized with constant speed, for any $0<\varepsilon<1$. The sequence $(\tilde{Q}_\varepsilon)_{0<\varepsilon<1}\subset\mathcal{P}(C([0,1],\Omega'))$ is tight: indeed, if we denote by
\begin{equation*}
    H([0,1],\Omega'):=\left\{\sigma:[0,1]\to\Omega':\sigma\textnormal{ is horizontal}\right\}
\end{equation*}
and we consider $K>0$, then the set $\left\{\sigma\in H([0,1],\Omega'):|\dot{\sigma}|_H\leq K\right\}$ is compact by the Ascoli-Arzelà Theorem. Moreover	
\begin{multline*}
	\tilde{Q}_\varepsilon\left(\left\{\sigma\in H([0,1],\Omega'):|\dot{\sigma}|_H>K\right\}\right)=\tilde{Q}_\varepsilon\left(\left\{\sigma\in H([0,1],\Omega'):\ell_{SR}(\sigma)>K\right\}\right)\\\leq \frac{1}{(1+\varepsilon)K}\int_{\Omega'}i_{Q_\varepsilon}(x)dx\leq\frac{c}{K}\|\textbf{w}\|_{L^p(\Omega,H\Omega)},
\end{multline*}
with $c>0$ independent of $\varepsilon$ thanks to \eqref{5dicembre}.
Therefore, $(\tilde{Q}_\varepsilon)_{0<\varepsilon<1}$ is tight, and thus admits a (not relabeled) subsequence that weakly converges to some $Q\in\mathcal{P}(C([0,1],\Omega'))$. It is obvious that $Q$ is concentrated on curves valued in $\overline{\Omega}$, therefore $Q\in \mathcal{P}(C([0,1],\overline{\Omega}))$. Moreover,  following the same proof of \cite[Lemma 2.8]{Santambrogio1}, one can prove that $Q$ is concentrated on $H$.
	
From the fact that  $\tilde{Q}_\varepsilon\rightharpoonup Q$, as $\varepsilon\to0$, and \eqref{approx1} we can conclude that
\begin{equation*}
	Q\in\mathcal{Q}_H(\mu,\nu).
\end{equation*}
In the end, $i_{\tilde{Q}_\varepsilon}=|\textbf{w}^\varepsilon|_H/(1+\varepsilon)$ converges to $|\textbf{w}|_H$ in $L^p(\Omega')$, hence it converges to the same limit in $\mathcal{M}_+(\Omega')$. Then, for any $\varphi\in C(\Omega')$, $\varphi\geq0$,
\begin{equation*}
\aligned	\int_{\Omega'}\varphi(x)|\textbf{w}(x)|_Hdx&=\lim_{\varepsilon\to0}\frac{1}{1+\varepsilon}\int_{\Omega'}\varphi(x)|\textbf{w}^\varepsilon(x)|_Hdx\\&=\lim_{\varepsilon\to0}\int_{C([0,1],\Omega')}\left(\int_0^1\varphi(\sigma(t))|\dot{\sigma}(t)|_Hdt\right)d\tilde{Q}_\varepsilon(\sigma)\\&\geq \int_{C([0,1],\Omega')}\left(\int_0^1\varphi(\sigma(t))|\dot{\sigma}(t)|_Hdt\right) dQ(\sigma)=\int_{\Omega'}\varphi(x) di_Q(x),
\endaligned
\end{equation*}
by the lower semicontinuity of the map $Q\mapsto \int\left(\int\varphi(\sigma(t))|\dot{\sigma}(t)|_Hdt\right) dQ(\sigma)$ with respect to the weak convergence of measures. In particular this means that  
\begin{equation*}
    i_Q\leq |\textbf{w}|_H\in L^p(\Omega),
\end{equation*}
hence 
\begin{equation*}
    Q\in\mathcal{Q}_H^p(\mu,\nu).
\end{equation*}
By using the monotonicity of $G$ 
\begin{equation*}
	\eqref{lepbme1chap4}\leq\int_\Omega G(x,i_Q(x))dx\leq \int_\Omega \mathcal{G}(x,\textbf{w}(x))dx=\eqref{generalbeckmann}.
\end{equation*}

\end{proof}

An alternative proof of the previous theorem can be found in the literature, specifically in \cite[Theorem 3.2]{Brasco}, and it is also applicable in $\mathbb{H}^n$. This proof relies on the well-known superposition principle and requires additional assumptions on the measures $\mu$ and $\nu$, namely that they are $p$-summable and bounded below. This approach does not require going through the regularization procedure of solutions to \eqref{generalbeckmann}. More precisely, let us suppose that $$\mu,\nu\in L^p(\Omega) \textnormal{ and bounded by below;}$$
take a minimizer $\textnormal{\textbf{w}}$ of \eqref{generalbeckmann} and consider the time-depending horizontal vector field
\begin{equation}
    \hat{\textnormal{\textbf{w}}}(t,x)=\frac{\textnormal{\textbf{w}}(x)}{(1-t)\mu(x)+t\nu(x)}.
\end{equation}
We know that $\mu-\nu=\textnormal{div}_H\textnormal{\textbf{w}}=\textnormal{div}_\epsilon\textnormal{\textbf{w}}$, for any $\epsilon>0$: hence, the linear interpolating curve $\overline\mu(t,x)=(1-t)\mu(x)+t\nu(x)$ is a positive measure-valued distributional solution of the continuity equation \eqref{11marzo} in $\mathbb{R}^{2n+1}$, with initial datum $\lambda(0,\cdot)=\mu$. Since $\textnormal{\textbf{w}}$ is a solution to \eqref{generalbeckmann}, then
\begin{equation*}
	\int_0^1\int_\Omega|\hat{\textnormal{\textbf{w}}}(t,x)|_\epsilon d\overline\mu(t,x)dt=\int_\Omega |{\textnormal{\textbf{w}}}(x)|_H dx<+\infty,
\end{equation*}
and we can apply \cite[Theorem 5.8]{Bernard}: it follows that there exists a probability measure $Q\in\mathcal{P}(C([0,1],\overline\Omega))$, such that $Q(AC([0,1],\overline\Omega))=1$, 
\begin{equation*}
    \overline\mu(t,\cdot)=(e_t)_\#Q
\end{equation*}
and $Q$-a.e. curve $\sigma$ is an integral curve of $\hat{\textnormal{\textbf{w}}}(t,x)$. Then $Q$-a.e. $\sigma$ is horizontal curve because the vector field $\textnormal{\textbf{w}}(t,x)$ is horizontal, and hence $Q\in\mathcal{Q}_H(\mu,\nu)$.

Following the same computation as $\eqref{16marzo}$, one can prove that
\begin{align*}
    \int_{\overline{\Omega}}\varphi(x)di_{Q}(x)=\int_\Omega\varphi(x)|{\textnormal{\textbf{w}}}(x)|_Hdx,
\end{align*}
for every $\varphi\in C(\overline\Omega)$. This implies that $i_{Q}=|{\textnormal{\textbf{w}}}|_H$ and thus $Q\in\mathcal{Q}^p_H(\mu,\nu)$ and it solves \eqref{lepbme1chap4}, concluding the proof.

Proposition \ref{21marzo251} and Theorem \ref{9dicembre2}, combined with the existence of solutions for \eqref{generalbeckmann} as established by Theorem \ref{existenceofbeckmin}, lead to the following two straightforward consequences.

\begin{cor}\label{11maggio1}
$\mathcal{Q}_H^p(\mu,\nu)\not=\emptyset$ if, and only if, $\mu-\nu\in \left(HW^{1,q}\right)'_\diamond\left(\Omega\right)$.
\end{cor}

\begin{cor}\label{24marzo25}
If $\mathcal{Q}_H^p(\mu,\nu)\not=\emptyset$, or equivalently $\mu-\nu\in \left(HW^{1,q}\right)'_\diamond\left(\Omega\right)$, then
\begin{enumerate}
    \item $\eqref{generalbeckmann}=\eqref{lepbme1chap4}$ and both $\eqref{generalbeckmann}$ and $\eqref{lepbme1chap4}$ admit solutions;
    \item for every $Q$ solving \eqref{lepbme1chap4}, $\textnormal{\textbf{w}}_Q$ solves \eqref{generalbeckmann}; 
    \item for every $\textnormal{\textbf{w}}$ solving \eqref{generalbeckmann}, there is $Q_{\textnormal{\textbf{w}}}$ that solves \eqref{lepbme1chap4}.
\end{enumerate}
\end{cor}

\begin{Remark}
If the function $i\mapsto G(x,i)$ is strictly increasing, for any $x\in\overline\Omega$, then for any $Q\in\mathcal{Q}_H^p(\mu,\nu)$, $Q$ solves \eqref{lepbme1chap4}$\Longleftrightarrow \textnormal{\textbf{w}}_Q$ solves \eqref{generalbeckmann} and $|\textnormal{\textbf{w}}_Q|_H=i_Q$.

Indeed, Corollary \ref{24marzo25}, Lemma \ref{7marzo24} and the monotonicity of $G$ imply
\begin{equation*}
	\eqref{lepbme1chap4}=\eqref{generalbeckmann}\leq\int_\Omega G(x,|\textnormal{\textbf{w}}_Q(x)|_H)dx\leq \int_\Omega G(x, i_Q(x))dx,\quad \forall Q\in\mathcal{Q}^p_H(\mu,\nu).
\end{equation*}
Therefore the result follows.

On the other hand, the convexity assumption on $i\mapsto G(x,i)$, for any $x\in\overline\Omega$, ensures the existence of minimizers for both problems \eqref{generalbeckmann} and \eqref{lepbme1chap4}. However, this convexity plays no role in the proof of the equivalence \eqref{generalbeckmann} = \eqref{lepbme1chap4}.
\end{Remark}

\section*{Acknowledgements}
The authors would like to thank the anonymous referee for the careful reading of the draft and for the valuable suggestions, which significantly contributed to the improvement of this paper.

M. C. is supported by the project PRIN 2022 F4F2LH - CUP J53D23003760006, \textit{Regularity
problems in sub-Riemannian structures}.

A.C. acknowledges partial support from grants PID2020-112881GB-I00, PID2021-125021NAI00 (Spanish Government) and 2021-SGR-00071 (Catalan Government).

\printbibliography

\end{document}